\documentstyle[12pt]{article}
\def\Bbb R{{\rm \bf R}}
\def\proclaim#1{\vskip2mm{\bf #1}\em}
\def\endproclaim{\em \vskip2mm}
\def\tag#1{\eqno(#1)}
\def\gathered{\begin{array}{c}}
\def\endgathered{\end{array}}
\def\text{\mbox}

\begin{document}

\title {On finding a cavity in a thermoelastic body
using a single displacement measurement over a finite time interval on the surface
of the body}
\author{Masaru IKEHATA\footnote{
Laboratory of Mathematics,
Graduate School of Engineering,
Hiroshima University,
Higashihiroshima 739-8527, Japan}}
\maketitle

\begin{abstract}
A mathematical formulation of an estimation problem of a cavity inside a three-dimensional thermoelastic body
using time domain data is considered.
The governing equation of the problem is given by a system of equations in the linear theory of thermoelasticity
which is a coupled system of the elastic wave and heat equations.
A new version of the enclosure method in the time domain which is originally developed for the classical wave equation is established.
For a comparison, the results in the decoupled case are also given.


\noindent
AMS: 74J25, 35Q79, 74F05, 35R30, 35L05, 35K05, 35B40

\end{abstract}


\section{Introduction}

The purpose of this paper is to pursue the possibility of the enclosure method \cite{I1, E00} itself
in inverse obstacle problems governed by several partial differential equations in {\it time domain}.
In particular, the enclosure method using a {\it single}
pair of the input and output over a {\it finite time interval} has been developed for several inverse obstacle problems
whose governing equations are given by scalar wave equations, heat equations and the Maxwell system,
see \cite{I4, IW00, IEO2, IEO3, ICA, Iwall, IK2, IK5, IMP, IMax, IMax2}.
In \cite{IE4} the author added a new idea to this time domain enclosure method for finding the geometry of an unknown obstacle from
a single point of the graph of the so-called {\it response operator} on the {\it outer surface} of the domain in which the obstacle is embedded.
The paper is focused on the explanation of the idea and so the governing equation of the wave therein is 
a classical wave equation.  

In this paper, we consider the new case when the governing equation is given by a coupled system of hyperbolic and parabolic 
equations.  How does this new enclosure method work for the system?

We restrict ourself to considering a classical system in the linear theory of thermoelasticity
since it is a typical coupled system having a physical meaning.
It is a coupled system of the elastic wave and heat equations.
We refer the reader to the article \cite{C} which gives us whole knowledge
about the system together with several references.

Now let us formulate the problem concretely. 
We denote by $\Omega$ and $D$ the reference body and an unknown cavity
embedded in $\Omega$, respectively.
It is assumed that $\Omega\setminus\overline D$ is homogeneous and isotropic. 
In this paper, for simplicity and considering the results in \cite{D}, we assume that $\Omega$ is given by a bounded domain with $C^{\infty}$-boundary
and $D$ a nonempty bounded open subset of $\Omega$ with $C^{\infty}$-boundary such that
$\Omega\setminus\overline D$ is connected.  We use the same symbol $\mbox{\boldmath $\nu$}$ to denote both the outer unit normal vectors of $\partial D$
and $\partial\Omega$. 

Let $0<T<\infty$.  Given $f=f(x,t)$ and $\mbox{\boldmath $G$}=\mbox{\boldmath $G$}(x,t)$ with
$(x,t)\in\partial\Omega\times\,]0,\,T[$ and  let
$\mbox{\boldmath $u$}=
\mbox{\boldmath $u$}_{f,\mbox{\boldmath $G$}}(x,t)$
and $\vartheta=\vartheta_{f, \mbox{\boldmath $G$}}(x,t)$ with $(x,t)\in\,(\Omega\setminus\overline D)\times\,]0,\,T[$ 
denote the solutions of the following initial boundary value problem
$$\displaystyle
\left\{
\begin{array}{ll}
\displaystyle
\rho\partial_t^2\mbox{\boldmath $u$}-\mu\Delta\mbox{\boldmath $u$}
-(\lambda+\mu)\nabla(\nabla\cdot\mbox{\boldmath $u$})-m\nabla\vartheta=\mbox{\boldmath $0$}
& \text{in}\,(\Omega\setminus\overline D)\times\,]0,\,T[,\\
\\
\displaystyle
c\partial_t\vartheta-k\Delta\vartheta-m\theta_0\nabla\cdot\partial_t\mbox{\boldmath $u$}=0
& \text{in}\,(\Omega\setminus\overline D)\times\,]0,\,T[,\\
\\
\displaystyle
\mbox{\boldmath $u$}(x,0)=\mbox{\boldmath $0$} & \text{in}\,\Omega\setminus\overline D,
\\
\\
\displaystyle
\partial_t\mbox{\boldmath $u$}(x,0)=\mbox{\boldmath $0$} & \text{in}\,\Omega\setminus\overline D,\\
\\
\displaystyle
\vartheta(x,0)=0 & \text{in}\,\Omega\setminus\overline D,\\
\\
\displaystyle
s(\mbox{\boldmath $u$}, \vartheta)\mbox{\boldmath $\nu$}
=\mbox{\boldmath $0$} & \text{on}\,\partial D\times\,]0,\,T[,\\
\\
\displaystyle
k\nabla\vartheta\cdot\mbox{\boldmath $\nu$}=0
& \text{on}\,\partial D\times\,]0,\,T[,\\
\\
\displaystyle
s(\mbox{\boldmath $u$},\vartheta)
\mbox{\boldmath $\nu$}=\mbox{\boldmath $G$} & \text{on}\,\partial\Omega\times\,]0,\,T[,
\\
\\
\displaystyle
-k\nabla\vartheta\cdot\mbox{\boldmath $\nu$}=f & \text{on}\,\partial\Omega\times\,]0,\,T[,
\end{array}
\right.
\tag {1.1}
$$
where
$$\
\displaystyle
s(\mbox{\boldmath $u$},\vartheta)
=2\mu\,\text{Sym}\,\nabla\mbox{\boldmath $u$}+
\lambda(\nabla\cdot\mbox{\boldmath $u$})I_3+m\,\vartheta I_3.
$$
The constants
$\theta_0$ and $m$ are the reference temperature and stress-temperature modulus of the body
$\Omega\setminus\overline D$, respectively;
$k$ the conductivity; $\lambda$ and $\mu$ are Lam\'e modulus and shear modulus, respectively;
$\rho$ and $c$ the density and specific heat.   
$\vartheta$ denotes the temperature difference of the absolute temperature from the reference temperature
$\theta_0$;
$\mbox{\boldmath $u$}$ and $s(\mbox{\boldmath $u$}, \vartheta)\mbox{\boldmath $\nu$}$ the displacement vector field
and the surface traction, respectively.

It is assumed that $\rho$, $c$, $\theta_0$ and $k$
are known positive constants; $m$, $\lambda$ and $\mu$ are known constants and satisfy $m\not=0$,
$\mu>0$ and $3\lambda+2\mu>0$.

Before describing the problem, we specify the class where the solution of (1.1) lives.
We employ a general result on the unique solvability and regularity of the initial boundary value problem for the coupled system
of the parabolic and hyperbolic system with {\it inhomogeneous Neumann-type boundary condition} established in \cite{D}.
It is based on the Hille-Yoshida theorem.
By applying Theorem 2.1 in \cite{D} to the present case, we see that
the initial boundary value problem (1.1) has a unique solution such that
$$
\left\{
\begin{array}{l}
\displaystyle
\mbox{\boldmath $u$}\in C^2([0,\,T];L^2(\Omega\setminus\overline D)^3)
\cap C^1([0,\,T];H^1(\Omega\setminus\overline D)^3)\cap C^0([0,\,T];H^2(\Omega\setminus\overline D)^3),\\
\\
\displaystyle
\partial_t^{1+l}\partial_x^{\alpha}\mbox{\boldmath $u$}
\in L^2([0,\,T];H^{-1/2}(\partial(\Omega\setminus\overline D))),\,l+\vert\alpha\vert\le 1,
\\
\\
\displaystyle
\vartheta\in
C^1([0,\,T];L^2(\Omega\setminus\overline D))
\cap C^0([0,\,T];H^2(\Omega\setminus\overline D)),
\\
\\
\displaystyle
\partial_t\vartheta\in L^2([0,\,T];H^1(\Omega\setminus\overline D))
\end{array}
\right.
$$
provided

$\bullet$ $f\in C^0([0,\,T];H^{1/2}(\partial\Omega))$ with
$\partial_tf\in L^2([0,\,T];H^{1/2}(\partial\Omega))$
and $f(0)=0$;

$\bullet$ $\mbox{\boldmath $G$}\in C^0([0,\,T];H^{1/2}(\partial\Omega)^3)$ with
$\partial_t\mbox{\boldmath $G$}\in L^2([0,\,T];H^{1/2}(\partial\Omega)^3)$
and $\mbox{\boldmath $G$}(0)=0$.

\noindent
We say that the pair $(f,\mbox{\boldmath $G$})$ is {\it admissible} if
the conditions listed above are satisfied.

Note that, this framework which is based on the Hille-Yoshida theorem corresponds to that of \cite{IE4} in which 
a combination of a standard lifting argument
and the results in \cite{I}
about the unique solvability of the classical wave equation with {\it homogeneous} Neumann boundary condition 
are employed.  See also \cite{S} for the hyperbolic systems with {\it inhomogeneous} Neumann
boundary condition.
For the classical results about the direct problem
with several {\it homogeneous} boundary conditions for the equations of theoremoelasticity 
which is based on a method due to Vi\v sik, see \cite{DA}.

\noindent

In this paper, we consider the following problem.

$\quad$

{\bf\noindent Problem.}  Fix a large $T$ (to be determined later) and a single set of the admissible 
pair  $(f, \mbox{\boldmath $G$})$ (to be specified later).
Assume that set $D$ is  unknown.  Extract information about the
location and shape of $D$ from the displacement field $\mbox{\boldmath $u$}_{f,\mbox{\boldmath $G$}}(x,t)$
and temperature difference $\vartheta_{f,\mbox{\boldmath $G$}}(x,t)$ given at
all $x\in\partial\Omega$ and $t\in\,]0,\,T[$.

$\quad$

For this problem we employ the idea in the most recent enclosure method
developed in \cite{IE4}.
The method introduces so-called indicator functions.

Let $B$ be an open ball centered at $p$ with radius $\eta$
satisfying $\overline B\cap\overline\Omega=\emptyset$.
We think that radius $\eta$ of $B$ is very small.
Let $\mbox{\boldmath $v$}$ and $\Theta$ satisfy
$$\displaystyle
\left\{
\begin{array}{ll}
\displaystyle
\rho\partial_t^2\mbox{\boldmath $v$}-\mu\Delta\mbox{\boldmath $v$}
-(\lambda+\mu)\nabla(\nabla\cdot\mbox{\boldmath $v$})-m\nabla\Theta=\mbox{\boldmath $0$}
& \text{in $\Bbb R^3\times\,]0,\,T[$,}\\
\\
\displaystyle
c\partial_t\Theta-k\Delta\Theta-m\theta_0\nabla\cdot\partial_t\mbox{\boldmath $v$}=0
& \text{in $\Bbb R^3\times\,]0,\,T[$,}
\end{array}
\right.
\tag {1.2}
$$
and
$$
\displaystyle
\left\{
\begin{array}{ll}
\displaystyle
\mbox{\boldmath $v$}(x,0)=\mbox{\boldmath $0$} & \text{in $\Bbb R^3$,}\\
\\
\displaystyle
\text{supp}\,\partial_t\mbox{\boldmath $v$}(\,\cdot\,,0)\cup
\text{supp}\,\Theta(\,\cdot\,,0)\subset\overline B. &
\end{array}
\right.
\tag {1.3}
$$
Note that, at this stage we do not specify the form of $\partial_t\mbox{\boldmath $v$}(\,\cdot\,,0)$
and $\Theta(\,\cdot\,,0)$.

The simplified version of the enclosure method employs
special $f$ and $\mbox{\boldmath $G$}$ in (1.1) as follows.

Set
$$\begin{array}{ll}
\displaystyle
\mbox{\boldmath $G$}(\mbox{\boldmath $v$},\Theta)=
s(\mbox{\boldmath $v$}, \Theta)\mbox{\boldmath $\nu$}
& \text{on $\partial\Omega\times\,]0,\,T[$}
\end{array}
\tag {1.4}
$$
and
$$\begin{array}{ll}
\displaystyle
f(\mbox{\boldmath $v$},\Theta)=-k\nabla\Theta\cdot\mbox{\boldmath $\nu$} & 
\text{on $\partial\Omega\times ]0,\,T[$.}
\end{array}
\tag {1.5}
$$
Note that both $\mbox{\boldmath $G$}(\mbox{\boldmath $v$},\Theta)$ 
and $f(\mbox{\boldmath $v$},\Theta)$ do not contain any {\it large} parameter.

We assume that the pair $(f,\mbox{\boldmath $G$})$ given by (1.4) and (1.5) is admissible.

Then, we introduce two indicator functions which play the central role in this paper.

{\bf\noindent Definition 1.1.}
Let $\mbox{\boldmath $u$}$ and $\vartheta$ solve (1.1) with $(\mbox{\boldmath $G$}, f)=(\mbox{\boldmath $G$}(\mbox{\boldmath $v$},\Theta), f(\mbox{\boldmath $v$},\Theta))$
given by (1.4) and (1.5), respectively.
Let $\tau>0$ and define
$$\begin{array}{ll}
\displaystyle
I^1(\tau;\mbox{\boldmath $v$},\Theta)=\int_{\partial\Omega}
s(\mbox{\boldmath $w$}_0,\Xi_0)\mbox{\boldmath $\nu$}
\cdot (\mbox{\boldmath $w$}-\mbox{\boldmath $w$}_0)\,dS
\end{array}
\tag {1.6}
$$
and
$$\begin{array}{lll}
\displaystyle
I^2(\tau;\mbox{\boldmath $v$},\Theta)=\int_{\partial\Omega}
k\frac{\partial\Xi_0}{\partial\nu}
 (\Xi-\Xi_0)\,dS,
\end{array}
\tag {1.7}
$$
where
$$
\left\{
\begin{array}{ll}
\displaystyle
\mbox{\boldmath $w$}(x)=\mbox{\boldmath $w$}(x,\tau)=\int_0^Te^{-\tau t}\mbox{\boldmath $u$}(x,t)dt,
&
x\in\Omega\setminus\overline D,\\
\\
\displaystyle
\mbox{\boldmath $w$}_0(x)=\mbox{\boldmath $w$}_0(x,\tau)=\int_0^Te^{-\tau t}\mbox{\boldmath $v$}(x,t)dt,
&
x\in\Bbb R^3,
\end{array}
\right.
\tag {1.8}
$$
and
$$
\left\{
\begin{array}{lll}
\displaystyle
\Xi(x)=\Xi(x,\tau)=\int_0^Te^{-\tau t}\vartheta(x,t)dt,
&
x\in\Omega\setminus\overline D,
\\
\\
\displaystyle
\Xi_0(x)=\Xi_0(x,\tau)=\int_0^Te^{-\tau t}\Theta(x,t)dt,
&
x\in\Bbb R^3.
\end{array}
\right.
\tag {1.9}
$$
Some remarks are in order.

$\bullet$  The indicator functions together with functions $\mbox{\boldmath $w$}$,
$\mbox{\boldmath $w$}_0$, $\Xi$ and $\Xi_0$ depend on $T$.
However, for simplicity of description, we omit to show their dependence on $T$
explicitly.

$\bullet$  The function $\mbox{\boldmath $w$}$ in (1.6) is the trace $\mbox{\boldmath $w$}\vert_{\partial\Omega}$
of $\mbox{\boldmath $w$}$
given by (1.8) 
onto $\partial\Omega$ and we have
$$\displaystyle
\mbox{\boldmath $w$}\vert_{\partial\Omega}
=\int_0^T e^{-\tau t}\mbox{\boldmath $u$}(\,\cdot\,,t)\vert_{\partial\Omega}\,dt.
$$
The integral on this right-hand is the integration for the functions with
the values in a Banach space (\cite{DL}).
The same remark works also for $\Xi$ in (1.7).
Therefore, these indicator functions can be computed from the response $\mbox{\boldmath $u$}$ 
and $\vartheta$ on $\partial\Omega$ over
time interval $]0,\,T[$ which are the solutions  of (1.1) with 
$\mbox{\boldmath $G$}=\mbox{\boldmath $G$}(\mbox{\boldmath $v$},\Theta)$
and $f=f(\mbox{\boldmath $v$},\Theta)$.

$\bullet$  Using the Lumer-Phillips theorem \cite{Y}, one can show that, given the initial data
$\mbox{\boldmath $v$}(x,0)=\mbox{\boldmath $v$}_0\in H^2(\Bbb R^3)^3$, $\partial_t\mbox{\boldmath $v$}(x,0)=
\mbox{\boldmath $v$}_1\in H^1(\Bbb R^3)^3$ and $\Theta(x,0)=\Theta_0\in H^2(\Bbb R^3)$, there exists
a unique pair of $\mbox{\boldmath $v$}\in C^2([0,\infty[, L^2(\Bbb R^3)^3)
\cap C^1([0,\infty[, H^1(\Bbb R^3)^3)\cap C^0([0,\infty[, H^2(\Bbb R^3)^3)$
and $\Theta\in C^1([0,\,\infty[, L^2(\Bbb R^3))\cap C^0([0,\,\infty[, H^2(\Bbb R^3))$ satisfying (1.2).
However, in this paper, we do not employ this general fact since the desired solutions $\mbox{\boldmath $v$}$ and $\Theta$ 
of (1.2) satisfying (1.3) have been constructed from those of decoupled equations.

Now we state the main results of this paper.

\subsection{Coupled case}

In this subsection $m$ in (1.1) is an arbitrary real number, and needless to say,
the special case $m=0$ is not excluded.

Let $\mbox{\boldmath $a$}$ be an arbitrary unit vector.
Let $\mbox{\boldmath $\Phi$}$ solve
$$\left\{\begin{array}{ll}
\displaystyle
\rho\partial_t^2\mbox{\boldmath $\Phi$}-\mu\Delta\mbox{\boldmath $\Phi$}=\mbox{\boldmath $0$}
&
\text{in $\Bbb R^3\times\,]0,\,T[$,}\\
\\
\displaystyle
\mbox{\boldmath $\Phi$}(x,0)=\mbox{\boldmath $0$} 
&
\text{in $\Bbb R^3$,}\\
\\
\displaystyle
\partial_t\mbox{\boldmath $\Phi$}(x,0)=(\eta-\vert x-p\vert)^2\chi_B(x)\mbox{\boldmath $a$}
&
\text{in $\Bbb R^3$,}
\end{array}
\right.
$$
where $\chi_B$ denotes the characteristic function of $B$.
Since the function $\Bbb R^3\ni x\longmapsto (\eta-\vert x-p\vert)^2\chi_B(x)\in\Bbb R$ belongs to 
$H^2(\Bbb R^3)$, it is known that one can construct such $\mbox{\boldmath $\Phi$}$ in the class
$$\displaystyle
C^2([0,\,T], H^1(\Bbb R^3)^3)
\cap C^1([0,\,T], H^2(\Bbb R^3)^3)
\cap C([0,\,T], H^3(\Bbb R^3)^3).
$$
Set
$$
\displaystyle
\mbox{\boldmath $v$}_s=\nabla\times\mbox{\boldmath $\Phi$}
\in
C^2([0,\,T], L^2(\Bbb R^3)^3)
\cap C^1([0,\,T], H^1(\Bbb R^3)^3)
\cap C([0,\,T], H^2(\Bbb R^3)^3).
\tag {1.10}
$$
We have $\nabla\cdot\mbox{\boldmath $v$}_s=0$.
We see that (1.2) and (1.3) are satisfied with the pair $(\mbox{\boldmath $v$},\Theta)
=(\mbox{\boldmath $v$}_s,0)$;
the pair $(f,\mbox{\boldmath $G$})=
(0,s(\mbox{\boldmath $v$}_s,0)\mbox{\boldmath $\nu$})$ given by (1.4) and (1.5)
is admissible.
Note that the pair $(\mbox{\boldmath $v$},\Theta)=(\mbox{\boldmath $v$}_s,0)$ 
is a special version of the Deresiewicz-Zorski solution of the system (1.2), see page 330 in \cite{C}.

In this case we have $\Xi_0=0$ in $\Bbb R^3$ and (1.7) gives
$I^2(\tau;\mbox{\boldmath $v$}_s,0)=0$ for all $\tau$.   
This means that one cannot obtain any information
about $D$ from the indicator function $\tau\longmapsto I^2(\tau;\mbox{\boldmath $v$}_s,0)$.
However, another indicator function $I^1(\tau;\mbox{\boldmath $v$}_s,0)$ 
has the following asymptotic behaviour.

\proclaim{\noindent Theorem 1.1.}
\noindent
(i)  Let $T$ satisfy
$$\displaystyle
T>\sqrt{\frac{\rho}{\mu}}\left(2\text{dist}\,(D,B)-\text{dist}\,(\Omega,B)\right).
\tag {1.11}
$$
Then, there exists a positive number $\tau_0$ such that, for all $\tau\ge\tau_0$
$I^1(\tau;\mbox{\boldmath $v$}_s,0)>0$ and we have
$$\displaystyle
\lim_{\tau\longrightarrow\infty}\frac{1}{\tau}\log I^1(\tau;\mbox{\boldmath $v$}_s,0)
= -2\sqrt{\frac{\rho}{\mu}}\text{dist}\,(D,B).
\tag {1.12}
$$

\noindent
(ii)  We have
$$\displaystyle
\lim_{\tau\longrightarrow\infty}e^{\tau T}I^1(\tau;\mbox{\boldmath $v$}_s,0)
=
\left\{\begin{array}{ll}
\displaystyle
\infty & 
\text{if $\displaystyle T>2\sqrt{\frac{\rho}{\mu}}\text{dist}\,(D,B)$,}
\\
\\
\displaystyle
0 & 
\text{if $\displaystyle T<2\sqrt{\frac{\rho}{\mu}}\text{dist}\,(D,B)$.}
\end{array}
\right.
\tag {1.13}
$$

\noindent
(iii)  If $\displaystyle T=2\sqrt{\frac{\rho}{\mu}}\text{dist}\,(D,B)$, then we have, as $\tau\longrightarrow\infty$
$$\displaystyle
e^{\tau T}I^1(\tau;\mbox{\boldmath $v$}_s,0)=O(\tau^4).
\tag {1.14}
$$

\endproclaim

Note that $\mbox{\boldmath $G$}$ depends on $\mbox{\boldmath $a$}$ 
and there is no restriction on the direction of $\mbox{\boldmath $a$}$ relative
to the unit normal $\mbox{\boldmath $\nu$}$ at the points on $\partial D$
which are nearest to the center point of ball $B$.

Since $T$ in (i) has the constraint (1.11), formula (1.12) does not automatically imply
the validity of (1.13) in  the case when $T<2\sqrt{\rho/\mu}\,\text{dist}\,(D,B)$.
Note also that constraint (1.11) is reasonable since as pointed out in \cite{IW00}
(see also \cite{IE4}) we have
$$
\displaystyle
2\text{dist}\,(D,B)-\text{dist}\,(\Omega,B)
\ge
\inf\{\vert x-y\vert+\vert y-z\vert\,\vert\,x\in\partial B, y\in\partial D,z\in\partial\Omega\}.
$$
The bound $O(\tau^4)$ on (1.14) is a rough estimation.  The point is: it is just at most algebraic.

Theorem 1.1 gives a solution to Problem in the case when $f=0$
and $\mbox{\boldmath $G$}=\mbox{\boldmath $G$}(\mbox{\boldmath $v$}_s,0)$,
where $\mbox{\boldmath $v$}_s$ is given by (1.10).
Needless to say, the indicator function $I^1(\tau;\mbox{\boldmath $v$}_s,0)$ contains information about thermal effect.
However, remarkably enough, the choice of the prescribed traction $\mbox{\boldmath $G$}=s(\mbox{\boldmath $v$}_s,0)\mbox{\boldmath $\nu$}$
and heat flux $f=0$ in (1.1) enables us to {\it ignore} thermal effects inside the body.
More precisely, with the help of the identity $I^2(\tau;\mbox{\boldmath $v$}_s,0)=0$ one can control the asymptotic behaviour
of the indicator function $I^1(\tau;\mbox{\boldmath $v$}_s,0)$ as $\tau\longrightarrow\infty$
in terms of $\mbox{\boldmath $v$}_s$ only over time interval $]0,\,T[$.  See Proposition 3.1 in Section 3 and formula (3.3) in the proof.

It should be pointed out that Theorem 1.1 does not give any solution
to the case when the surface traction on the outer boundary
vanishes, that is, $\mbox{\boldmath $G$}=\mbox{\boldmath $0$}$.
This case shall be more difficult since it seems that system (1.2) does not posses a solution $(\mbox{\boldmath $v$},\Theta)$
in such a way that $G(\mbox{\boldmath $v$},\Theta)=\mbox{\boldmath $0$}$ on $\partial\Omega\times\,]0,\,T[$.

\subsection{Decoupled case}

The enclosure method presented in this paper covers also the decoupled case $m=0$.
In this subsection, for a comparison with the main result we present two results in that case.  

Let $m=0$.  
Let $\phi$ be the solution of
$$\left\{
\begin{array}{ll}
\rho\partial_t^2\phi-(\lambda+2\mu)\Delta\phi=0
&
\text{in $\Bbb R^3\times\,]0,\,T[$,}
\\
\\
\displaystyle
\phi(x,0)=0 & 
\text{in $\Bbb R^3$,}
\\
\\
\displaystyle
\partial_t\phi(x,0)=(\eta-\vert x-p\vert)^2\chi_B(x)
&
\text{in $\Bbb R^3$.}
\end{array}
\right.
$$
The class where $\phi$ belongs to is the same as all the components of $\mbox{\boldmath $\Phi$}$
in the preceding section.
Set
$$\displaystyle
\mbox{\boldmath $v$}_p=\nabla\phi.
\tag {1.15}
$$
The same comment for the class where $\mbox{\boldmath $v$}_p$ works also
and we have $\nabla\times\mbox{\boldmath $v$}_p=\mbox{\boldmath $0$}$.
We see that the pair $(\mbox{\boldmath $v$}, \Theta)= (\mbox{\boldmath $v$}_p,0)$ satisfies
(1.2) and (1.3) under the assumption $m=0$;
the pair $(f,\mbox{\boldmath $G$})=(0, s(\mbox{\boldmath $v$}_p,0)\mbox{\boldmath $\nu$})$ given by (1.4) and (1.5) is admissible.

\proclaim{\noindent Theorem 1.2.}  Let $m=0$.
\noindent
(i)  Let $T$ satisfy
$$\displaystyle
T>\sqrt{\frac{\rho}{\lambda+2\mu}}\left(2\text{dist}\,(D,B)-\text{dist}\,(\Omega,B)\right).
$$
Then, there exists a positive number $\tau_0$ such that, for all $\tau\ge\tau_0$
$I^1(\tau;\mbox{\boldmath $v$}_p,0)>0$ and we have
$$\displaystyle
\lim_{\tau\longrightarrow\infty}\frac{1}{\tau}\log I^1(\tau;\mbox{\boldmath $v$}_p,0)
= -2\sqrt{\frac{\rho}{\lambda+2\mu}}\text{dist}\,(D,B).
$$

\noindent
(ii)  We have
$$\displaystyle
\lim_{\tau\longrightarrow\infty}e^{\tau T}I^1(\tau;\mbox{\boldmath $v$}_p,0)
=
\left\{\begin{array}{ll}
\displaystyle
\infty & 
\text{if $\displaystyle T>2\sqrt{\frac{\rho}{\lambda+2\mu}}\text{dist}\,(D,B)$,}
\\
\\
\displaystyle
0 & 
\text{if $\displaystyle T<2\sqrt{\frac{\rho}{\lambda+2\mu}}\text{dist}\,(D,B)$.}
\end{array}
\right.
$$

\noindent
(iii)  If $\displaystyle T=2\sqrt{\frac{\rho}{\lambda+2\mu}}
\text{dist}\,(D,B)$, then we have, as $\tau\longrightarrow\infty$
$$\displaystyle
e^{\tau T}I^1(\tau;\mbox{\boldmath $v$}_p,0)=O(\tau^4).
$$

\endproclaim

Note that Theorem 1.1 covers also the case when $m=0$.  Thus, we have two results
in that case.

Let $\Theta_0$ solve
$$\left\{
\begin{array}{ll}
\displaystyle
c\partial_t\Theta-k\Delta\Theta=0
& \text{in $\Bbb R^3\times\,]0,\,T[$,}\\
\\
\displaystyle
\Theta(x,0)=(\eta-\vert x-p\vert)^2\chi_B(x) & \text{in $\Bbb R^3$.}
\end{array}
\right.
\tag {1.16}
$$
We see that the pair $(\mbox{\boldmath $v$},\Theta)=(\mbox{\boldmath $0$},\Theta_0)$
satisfies (1.2) and (1.3) under the assumption $m=0$;
the pair $(f,\mbox{\boldmath $G$})=(-k\nabla\Theta_0\cdot\mbox{\boldmath $\nu$},\mbox{\boldmath $0$})$
given by (1.4) and (1.5)
is admissible.

\proclaim{\noindent Theorem 1.3.}  Let $m=0$.

\noindent
(i)  Let $T$ be an arbitrary fixed positive number.
Then,
there exists a positive number $\tau_0$ such that, for all $\tau\ge\tau_0$
$I^2(\tau;\mbox{\boldmath $0$},\Theta_0)>0$ and we have
$$\displaystyle
\lim_{\tau\longrightarrow\infty}\frac{1}{\sqrt{\tau}}\log I^2(\tau;\mbox{\boldmath $0$},\Theta_0)
= -2\sqrt{\frac{c}{k}}\text{dist}\,(D,B).
$$

\noindent
(ii)  We have
$$\displaystyle
\lim_{\tau\longrightarrow\infty}e^{\sqrt{\tau} \,T}I^2(\tau;\mbox{\boldmath $0$},\Theta_0)
=
\left\{\begin{array}{ll}
\displaystyle
\infty & 
\text{if $\displaystyle T>2\sqrt{\frac{c}{k}}\text{dist}\,(D,B)$,}
\\
\\
\displaystyle
0 & 
\text{if $\displaystyle T<2\sqrt{\frac{c}{k}}\text{dist}\,(D,B)$.}
\end{array}
\right.
$$

\noindent
(iii)  If $\displaystyle T=2\sqrt{\frac{c}{k}}\text{dist}\,(D,B)$, then we have, as $\tau\longrightarrow\infty$
$$\displaystyle
e^{\sqrt{\tau}\,T}I^2(\tau;\mbox{\boldmath $0$},\Theta_0)=O(\tau^2).
$$

\endproclaim

Note that, in (i) of Theorem 1.3 there is no restriction on the size of $T$
same as several results in \cite{I4, IFR, IK1, IK2, IK5, II}.
This suggests the {\it infinite propagation} speed of the signal governed by the heat equation
indirectly.  By virtue of arbitrariness of $T$ one can easily deduce (ii) from (i).  
However, (i) in Theorems 1.1 and 1.2
there are restrictions on the size of $T$.  Thus to derive the later half of (ii) in Theorems 1.1 and 1.2
is independent of (i) in those theorems.

Since the proof of Theorems 1.2-1.3 is easier than that of Theorem 1.1 and can be done 
as Theorem 1.1 in \cite{IE4}, we omit to describe their proof.

\section{Preliminaries}

\subsection{Decomposition formulae of the indicator functions}

In this and next subsection section, 
the form of the pair $(\mbox{\boldmath $v$},\Theta)$ satisfying (1.2) and (1.3)
is not specified.

Set
$$
\begin{array}{ll}
\displaystyle
\mbox{\boldmath $R$}=\mbox{\boldmath $w$}-\mbox{\boldmath $w$}_0,
&
\displaystyle
\Sigma=\Xi-\Xi_0,
\end{array}
$$
where $\mbox{\boldmath $w$}$, $\Xi$, $\mbox{\boldmath $w$}_0$ and
$\Xi_0$ are given by (1.8) and (1.9).

It follows from (1.1) that $\mbox{\boldmath $w$}$ and $\Xi$ satisfies
$$\left\{
\begin{array}{ll}
\displaystyle
\mu\Delta\mbox{\boldmath $w$}+(\lambda+\mu)\nabla(\nabla\cdot\mbox{\boldmath $w$})
-\rho\tau^2\mbox{\boldmath $w$}
+m\nabla\Xi=\rho e^{-\tau T}\mbox{\boldmath $F$}& \text{in}\,\Omega\setminus\overline D,\\
\\
\displaystyle
k\Delta\Xi-c\tau\Xi+m\theta_0\tau\nabla\cdot\mbox{\boldmath $w$}
=e^{-\tau T}h & \text{in}\,\Omega\setminus\overline D,
\\
\\
\displaystyle
s(\mbox{\boldmath $w$},\Xi)\mbox{\boldmath $\nu$}
=s(\mbox{\boldmath $w$}_0,\Xi_0)\mbox{\boldmath $\nu$} & \text{on}\,\partial\Omega,\\
\\
\displaystyle
-k\nabla\Xi\cdot\mbox{\boldmath $\nu$}=
-k\nabla\Xi_0\cdot\mbox{\boldmath $\nu$}
& \text{on}\,\partial\Omega,
\\
\\
\displaystyle
s(\mbox{\boldmath $w$},\Xi)\mbox{\boldmath $\nu$}
=\mbox{\boldmath $0$} & \text{on}\,\partial D,
\\
\\
\displaystyle
-k\nabla\Xi\cdot\mbox{\boldmath $\nu$}=0 &
\text{on}\,\partial D,
\end{array}
\right.
\tag {2.1}
$$
where
$$
\left\{
\begin{array}{ll}
\displaystyle
\mbox{\boldmath $F$}=\mbox{\boldmath $F$}(x,\tau)=\partial_t\mbox{\boldmath $u$}(x,T)
+\tau\mbox{\boldmath $u$}(x,T), & x\in\Omega\setminus\overline D,\\
\\
\displaystyle
h=c\vartheta(x,T)-m\theta_0\nabla\cdot\mbox{\boldmath $u$}(x,T),
&
 x\in\Omega\setminus\overline D.
\end{array}
\right.
\tag {2.2}
$$
It follows from (1.2) that the $\mbox{\boldmath $w$}_0$ and $\Theta_0$
satisfy
$$\left\{
\begin{array}{ll}
\displaystyle
\mu\Delta\mbox{\boldmath $w$}_0+(\lambda+\mu)\nabla(\nabla\cdot\mbox{\boldmath $w$}_0)
-\rho\tau^2\mbox{\boldmath $w$}_0
+m\nabla\Xi_0+\rho\mbox{\boldmath $v$}_0
=\rho e^{-\tau T}\mbox{\boldmath $F$}_0& \text{in}\,\Bbb R^3,\\
\\
\displaystyle
k\Delta\Xi_0-c\tau\Xi_0+m\theta_0\tau\nabla\cdot\mbox{\boldmath $w$}_0+cf_0
=e^{-\tau T}h_0& \text{in}\,\Bbb R^3,
\end{array}
\right.
\tag {2.3}
$$
where
$$
\left\{
\begin{array}{ll}
\displaystyle
\mbox{\boldmath $v$}_0(x)=\partial_t\mbox{\boldmath $v$}(x,0),
& x\in\Bbb R^3,\\
\\
\displaystyle
f_0(x)=\Theta(x,0),
& x\in\Bbb R^3,
\\
\\
\displaystyle
\mbox{\boldmath $F$}_0=\mbox{\boldmath $F$}_0(x,\tau)=\partial_t\mbox{\boldmath $v$}(x,T)
+\tau\mbox{\boldmath $v$}(x,T), & x\in\Bbb R^3,\\
\\
\displaystyle
h_0=c\Theta(x,T)-m\theta_0\nabla\cdot\mbox{\boldmath $v$}(x,T),
&
 x\in\Bbb R^3.
\end{array}
\right.
\tag {2.4}
$$
Recall that the supports of $\mbox{\boldmath $v$}_0$ and $f_0$ are contained in $\overline B$
and $B$ satisfies $\overline B\cap\overline\Omega=\emptyset$.  See assumption (1.3).

Then, integration by parts together with the first equations on (2.1) and (2.3) yields
$$\begin{array}{l}
\,\,\,\,\,\,
\displaystyle
\int_{\partial\Omega}
\left(s(\mbox{\boldmath $w$}_0,\Xi_0)\mbox{\boldmath $\nu$}\cdot
\mbox{\boldmath $w$}-
s(\mbox{\boldmath $w$},\Xi)\mbox{\boldmath $\nu$}\cdot
\mbox{\boldmath $w$}_0\right)\,dS\\
\\
\displaystyle
=\int_{\partial D}s(\mbox{\boldmath $w$}_0,\Xi_0)\mbox{\boldmath $\nu$}\cdot
\mbox{\boldmath $w$}\,dS
+m\int_{\Omega\setminus\overline D}(\Xi_0\nabla\cdot\mbox{\boldmath $w$}
-\Xi\nabla\cdot\mbox{\boldmath $w$}_0)\,dx\\
\\
\,\,\,
\displaystyle
+e^{-\tau T}\rho
\int_{\Omega\setminus\overline D}
\left
(\mbox{\boldmath $F$}_0\cdot\mbox{\boldmath $w$}-\mbox{\boldmath $F$}\cdot\mbox{\boldmath $w$}_0
\right)
dx.
\end{array}
$$
Hence
$$\begin{array}{ll}
\displaystyle
I^1(\tau;\mbox{\boldmath $v$},\Theta)
&
\displaystyle
=\int_{\partial D}s(\mbox{\boldmath $w$}_0,\Xi_0)\mbox{\boldmath $\nu$}\cdot
\mbox{\boldmath $w$}\,dS
+m\int_{\Omega\setminus\overline D}(
-\Sigma\nabla\cdot\mbox{\boldmath $w$}_0
+
\Xi_0\nabla\cdot\mbox{\boldmath $R$})\,dx
\\
\\
\displaystyle
&
\displaystyle
\,\,\,
+e^{-\tau T}\rho
\int_{\Omega\setminus\overline D}
\left
(\mbox{\boldmath $F$}_0\cdot\mbox{\boldmath $w$}-\mbox{\boldmath $F$}\cdot\mbox{\boldmath $w$}_0
\right)
dx.
\end{array}
\tag {2.5}
$$
This is the first representation of indicator function $I^1(\tau;\mbox{\boldmath $v$},\Theta)$.

It follows from the second equations on (2.1) and (2.3) that
$$
\begin{array}{l}
\,\,\,\,\,\,
\displaystyle
\int_{\partial\Omega}
\left(k\frac{\partial\Xi_0}{\partial\nu}\Xi-k\frac{\partial\Xi}{\partial\nu}\Xi_0\right)\,dS\\
\\
\displaystyle
=
\int_{\partial D}k\frac{\partial\Xi_0}{\partial\nu}\Xi\,dS
+m\theta_0\tau
\int_{\Omega\setminus\overline D}
\left(\Xi_0\nabla\cdot\mbox{\boldmath $w$}-\Xi\nabla\cdot\mbox{\boldmath $w$}_0\right)\,dx\\
\\
\displaystyle
\,\,\,
+
e^{-\tau T}
\int_{\Omega\setminus\overline D}
(h_0\Xi-h\Xi_0)\,dx.
\end{array}
$$
This yields
$$
\begin{array}{ll}
\displaystyle
I^2(\tau;\mbox{\boldmath $v$},\Theta)
&
\displaystyle
=\int_{\partial D}k\frac{\partial\Xi_0}{\partial\nu}\Xi\,dS
+m\theta_0\tau
\int_{\Omega\setminus\overline D}
\left(-\Sigma\nabla\cdot\mbox{\boldmath $w$}_0+\Xi_0\nabla\cdot\mbox{\boldmath $R$}\right)\,dx
\\
\\
\displaystyle
&
\displaystyle
\,\,\,
+
e^{-\tau T}
\int_{\Omega\setminus\overline D}
(h_0\Xi-h\Xi_0)\,dx.
\end{array}
\tag {2.6}
$$
This is the first representation of indicator function $I^2(\tau;\mbox{\boldmath $v$},\Theta)$.

Next we decompose the first term on the right-hand side of
(2.5).  The result yields the following decomposition formula for 
$I^1(\tau;\mbox{\boldmath $v$},\Theta)$.

\proclaim{\noindent Proposition 2.1.}
We have
$$\begin{array}{ll}
\displaystyle
I^1(\tau;\mbox{\boldmath $v$},\Theta)
&
\displaystyle
=J(\tau)+m\int_D\Xi_0\nabla\cdot\mbox{\boldmath $w$}_0
+
E(\tau)
+m\int_{\Omega\setminus\overline D}\Sigma\nabla\cdot\mbox{\boldmath $R$}\,dx\\
\\
\displaystyle
&
\displaystyle
\,\,\,
+
m\int_{\Omega\setminus\overline D}(
-\Sigma\nabla\cdot\mbox{\boldmath $w$}_0
+
\Xi_0\nabla\cdot\mbox{\boldmath $R$})\,dx
+{\cal R}^1(\tau),
\end{array}
\tag {2.7}
$$
where
$$\displaystyle
J(\tau)=
\int_D\left(2\mu\left\vert\text{Sym}\,\nabla\mbox{\boldmath $w$}_0\right\vert^2
+\lambda\vert\nabla\cdot\mbox{\boldmath $w$}_0\vert^2+\rho\tau^2\vert\mbox{\boldmath $w$}_0\vert^2\right)\,dx,
\tag {2.8}
$$
$$\displaystyle
E(\tau)=\int_{\Omega\setminus\overline D}\left(2\mu\left\vert\text{Sym}\,\nabla\mbox{\boldmath $R$}\right\vert^2
+\lambda\vert\nabla\cdot\mbox{\boldmath $R$}\vert^2+\rho\tau^2\vert\mbox{\boldmath $R$}\vert^2\right)\,dx
\tag {2.9}
$$
and
$$
\displaystyle
{\cal R}^1(\tau)
=\rho e^{-\tau T}\left\{
\int_D\mbox{\boldmath $F$}_0\cdot\mbox{\boldmath $w$}_0\,dx
+\int_{\Omega\setminus\overline D}\mbox{\boldmath $F$}\cdot\mbox{\boldmath $R$}\,dx
+\int_{\Omega\setminus\overline D}
(\mbox{\boldmath $F$}_0-\mbox{\boldmath $F$})\cdot\mbox{\boldmath $w$}_0
dx\right\}.
\tag {2.10}
$$

\endproclaim

{\it\noindent Proof.}
Since $\overline B\cap\overline\Omega=\emptyset$, 
from the first equations on (2.1) and (2.3)  we see that
the $\mbox{\boldmath $R$}$ satisfies
$$\left\{
\begin{array}{ll}
\displaystyle
\mu\Delta\mbox{\boldmath $R$}+(\lambda+\mu)\nabla(\nabla\cdot\mbox{\boldmath $R$})
-\rho\tau^2\mbox{\boldmath $R$}
+m\nabla\Sigma=\rho e^{-\tau T}(\mbox{\boldmath $F$}-\mbox{\boldmath $F$}_0) & \text{in}\,\Omega\setminus\overline D,\\
\\
\displaystyle
s(\mbox{\boldmath $R$},\Sigma)\mbox{\boldmath $\nu$}
=\mbox{\boldmath $0$} & \text{on}\,\partial\Omega,
\\
\\
\displaystyle
\displaystyle
s(\mbox{\boldmath $R$},\Sigma)\mbox{\boldmath $\nu$}=
-s(\mbox{\boldmath $w$}_0,\Xi_0)\mbox{\boldmath $\nu$}
 & \text{on}\,\partial D.
\end{array}
\right.
\tag {2.11}
$$
Then, one can write
$$\displaystyle
\int_{\partial D}s(\mbox{\boldmath $w$}_0,\Xi_0)\mbox{\boldmath $\nu$}\cdot
\mbox{\boldmath $w$}\,dS
=\int_{\partial D}s(\mbox{\boldmath $w$}_0,\Xi_0)\mbox{\boldmath $\nu$}\cdot
\mbox{\boldmath $w$}_0\,dS
-\int_{\partial D}s(\mbox{\boldmath $R$},\Sigma)\mbox{\boldmath $\nu$}\cdot
\mbox{\boldmath $R$}\,dS.
$$
It follows from the first equation on (2.3) that
$$
\displaystyle
\int_{\partial D}s(\mbox{\boldmath $w$}_0,\Xi_0)\mbox{\boldmath $\nu$}\cdot
\mbox{\boldmath $w$}_0\,dS
=J(\tau)+m\int_D\Xi_0\nabla\cdot\mbox{\boldmath $w$}_0\,dx
+\rho e^{-\tau T}\int_D\mbox{\boldmath $F$}_0\cdot\mbox{\boldmath $w$}_0\,dx.
$$
It follows from (2.11) that
$$\begin{array}{ll}
\displaystyle
-\int_{\partial D}s(\mbox{\boldmath $R$},\Sigma)\mbox{\boldmath $\nu$}\cdot
\mbox{\boldmath $R$}\,dS
&
\displaystyle
=\int_{\partial(\Omega\setminus\overline D)}s(\mbox{\boldmath $R$},\Sigma)\mbox{\boldmath $\nu$}\cdot
\mbox{\boldmath $R$}\,dS\\
\\
\displaystyle
&
\displaystyle
=E(\tau)+m\int_{\Omega\setminus\overline D}\Sigma\nabla\cdot\mbox{\boldmath $R$}\,dx
+\rho e^{-\tau T}\int_{\Omega\setminus\overline D}(\mbox{\boldmath $F$}
-
\mbox{\boldmath $F$}_0)\cdot\mbox{\boldmath $R$}\,dx.
\end{array}
$$
Thus we obtain
$$\begin{array}{ll}
\displaystyle
\int_{\partial D}s(\mbox{\boldmath $w$}_0,\Xi_0)\mbox{\boldmath $\nu$}\cdot
\mbox{\boldmath $w$}\,dS
&
\displaystyle
=J(\tau)+m\int_D\Xi_0\nabla\cdot\mbox{\boldmath $w$}_0\,dx
+E(\tau)+m\int_{\Omega\setminus\overline D}\Sigma\nabla\cdot\mbox{\boldmath $R$}\,dx
\\
\\
\displaystyle
\,\,\,
&
\,\,\,
\displaystyle
+\rho e^{-\tau T}
\left\{
\int_D\mbox{\boldmath $F$}_0\cdot\mbox{\boldmath $w$}_0\,dx+
\int_{\Omega\setminus\overline D}(\mbox{\boldmath $F$}
-
\mbox{\boldmath $F$}_0)\cdot\mbox{\boldmath $R$}\,dx\right\}.
\end{array}
$$
Then a combination of this and (2.5) gives (2.7).

\noindent
$\Box$

{\bf\noindent Remark 2.1.}
We have, for all real $3\times 3$-matrix $A$ 
$$
\displaystyle
2\mu\left\vert\text{Sym}\,A-\frac{\text{trace}\,A}{3}I_3\right\vert^2+
\frac{3\lambda+2\mu}{3}\vert\text{trace}\,A\vert^2
=2\mu\vert\text{Sym}\,A\vert^2+\lambda\vert\text{trace}\,A\vert^2.
\tag {2.12}
$$
See page 85 in \cite{Gu}.  
Thus, both $J(\tau)$ and $E(\tau)$ given by (2.8) and (2.9), respectively are nonnegative under the assumption
$\mu>0$ and $3\lambda+2\mu>0$.

\proclaim{\noindent Proposition 2.2.}
We have

$$
\begin{array}{ll}
\displaystyle
I^2(\tau;\mbox{\boldmath $v$},\Theta)
&
\displaystyle
=j(\tau)-m\theta_0\tau\int_D\Xi_0\nabla\cdot\mbox{\boldmath $w$}_0\,dx
+e(\tau)-m\theta_0\tau\int_{\Omega\setminus\overline D}\Sigma\nabla\cdot\mbox{\boldmath $R$}\,dx
\\
\\
\displaystyle
\,\,\,
&
\displaystyle
\,\,\,
+m\theta_0\tau
\int_{\Omega\setminus\overline D}
\left(-\Sigma\nabla\cdot\mbox{\boldmath $w$}_0+\Xi_0\nabla\cdot\mbox{\boldmath $R$}\right)\,dx
+{\cal R}^2(\tau),
\end{array}
\tag {2.13}
$$
where
$$\displaystyle
j(\tau)=
\int_D
\left(
k\vert\nabla\Xi_0\vert^2
+c\tau\vert\Xi_0\vert^2\right)\,dx,
\tag {2.14}
$$
$$\displaystyle
e(\tau)=\int_{\Omega\setminus\overline D}\left(
k\vert\nabla\Sigma\vert^2
+c\tau\vert\Sigma\vert^2\right)\,dx
\tag {2.15}
$$
and
$$
\displaystyle
{\cal R}^2(\tau)
=e^{-\tau T}\left\{
\int_Dh_0\Xi_0\,dx
+\int_{\Omega\setminus\overline D}h\Sigma\,dx
+\int_{\Omega\setminus\overline D}
(h_0-h)\Xi_0
dx\right\}
$$

\endproclaim

{\it\noindent Proof.}
Since $\overline B\cap\overline\Omega=\emptyset$,
it follows from some of equations on (2.1) and (2.3) that
the $\Sigma$ satisfies
$$\left\{
\begin{array}{ll}
\displaystyle
k\Delta\Sigma-c\tau\Sigma+m\theta_0\tau\nabla\cdot\mbox{\boldmath $R$}
=e^{-\tau T}(h-h_0) & \text{in}\,\Omega\setminus\overline D,\\
\\
\displaystyle
-k\nabla\Sigma\cdot\mbox{\boldmath $\nu$}=0 & \text{on}\,\partial\Omega,\\
\\
\displaystyle
-k\nabla\Sigma\cdot\mbox{\boldmath $\nu$}
=k\nabla\Xi_0\cdot\mbox{\boldmath $\nu$}
 & \text{on}\,\partial D.
\end{array}
\right.
\tag {2.16}
$$
Then, one can write
$$\displaystyle
\int_{\partial D}k\frac{\partial\Xi_0}{\partial\mbox{\boldmath $\nu$}}\Xi\,dS
=\int_{\partial D}k\frac{\partial\Xi_0}{\partial\mbox{\boldmath $\nu$}}\Xi_0\,dS
-\int_{\partial D}k\frac{\partial\Sigma}{\partial\mbox{\boldmath $\nu$}}\Sigma\,dS.
$$
It follows from the second equation on  (2.3) that
$$
\displaystyle
\int_{\partial D}k\frac{\partial\Xi_0}{\partial\mbox{\boldmath $\nu$}}\Xi_0\,dS
=j(\tau)
-m\theta_0\tau\int_D\Xi_0\nabla\cdot\mbox{\boldmath $w$}_0\,dx
+e^{-\tau T}\int_Dh_0\Xi_0\,dx.
$$
It follows from (2.16) that
$$\begin{array}{ll}
\displaystyle
-\int_{\partial D}k\frac{\partial\Sigma}{\partial\mbox{\boldmath $\nu$}}\Sigma\,dS
&
\displaystyle
=
\int_{\partial\,(\Omega\setminus\overline D)}
k\frac{\partial\Sigma}{\partial\mbox{\boldmath $\nu$}}\Sigma\,dS\\
\\
\displaystyle
&
\displaystyle
=e(\tau)
-m\theta_0\tau\int_{\Omega\setminus\overline D}\Sigma\nabla\cdot\mbox{\boldmath $R$}\,dx
+e^{-\tau T}\int_{\Omega\setminus\overline D}(h-h_0)\Sigma\,dx.
\end{array}
$$
Thus we obtain
$$\begin{array}{ll}
\displaystyle
\int_{\partial D}k\frac{\partial\Xi_0}{\partial\mbox{\boldmath $\nu$}}\Xi\,dS
&
\displaystyle
=j(\tau)-m\theta_0\tau\int_D\Xi_0\nabla\cdot\mbox{\boldmath $w$}_0\,dx
+e(\tau)-m\theta_0\tau\int_{\Omega\setminus\overline D}\Sigma\nabla\cdot\mbox{\boldmath $R$}\,dx
\\
\\
\displaystyle
\,\,\,
&
\,\,\,
\displaystyle
+e^{-\tau T}\left\{\int_Dh_0\Xi_0dx+\int_{\Omega\setminus\overline D}(h-h_0)\mbox{\boldmath $\Sigma$}\,dx\right\}.
\end{array}
$$
Then a combination of this and (2.6) gives (2.13).

\noindent
$\Box$

As a direct consequence of Propositions 2.1 and 2.2 we obtain
$$\begin{array}{l}
\displaystyle
\,\,\,\,\,\,
\theta_0\tau I^1(\tau;\mbox{\boldmath $v$},\Theta)
+I^2(\tau;\mbox{\boldmath $v$},\Theta)\\
\\
\displaystyle
=(j(\tau)+\theta_0\tau\,J(\tau))+(e(\tau)+\theta_0\tau\,E(\tau))
\\
\\
\displaystyle
\,\,\,
+m\theta_0\tau\,\int_{\Omega\setminus\overline D}
(-\Sigma\nabla\cdot\mbox{\boldmath $w$}_0+\Xi_0\nabla\cdot\mbox{\boldmath $R$})\,dx
+\theta_0\tau\,{\cal R}(\tau),
\end{array}
\tag {2.17}
$$
where
$$\displaystyle
{\cal R}(\tau)={\cal R}^1(\tau)+\frac{1}{\theta_0\tau}\,{\cal R}^2(\tau).
\tag {2.18}
$$
Note that the first and second terms in the right-hand side on (2.17) are non negative
and the third term contains $\nabla\cdot\mbox{\boldmath $w$}_0$ and $\Xi_0$ linearly.
In short, we choose special $(\mbox{\boldmath $v$},\Theta)$ in Theorem 1.1 in such a way that
this third term vanishes. 

In the next subsection we give an upper bound on $e(\tau)+\theta_0\tau\,E(\tau)$ in terms of $j(\tau)$ and $J(\tau)$.

\subsection{Basic estimate}

\proclaim{\noindent Proposition 2.3.}
We have, as $\tau\longrightarrow\infty$
$$\displaystyle
e(\tau)+\theta_0\tau E(\tau)
=O(\tau j(\tau)+\tau^3 J(\tau)+\tau^3 e^{-2\tau T}).
\tag {2.19}
$$

\endproclaim

{\it\noindent Proof.}
It follows from (2.16) that
$$\begin{array}{ll}
\displaystyle
\int_{\partial D}k\frac{\partial\Xi_0}{\partial\mbox{\boldmath $\nu$}}\Sigma\,dS
&
\displaystyle
=e(\tau)
-m\theta_0\tau\int_{\Omega\setminus\overline D}\Sigma\nabla\cdot\mbox{\boldmath $R$}\,dx
+e^{-\tau T}\int_{\Omega\setminus\overline D}(h-h_0)\Sigma\,dx.
\end{array}
\tag {2.20}
$$
On the other hand, it follows from (2.11) that
$$\begin{array}{ll}
\displaystyle
\theta_0\tau\int_{\partial D}s(\mbox{\boldmath $w$}_0,\Xi_0)\mbox{\boldmath $\nu$}\cdot
\mbox{\boldmath $R$}\,dS
&
\displaystyle
=\theta_0\tau E(\tau)+m\theta_0\tau\int_{\Omega\setminus\overline D}\Sigma\nabla\cdot\mbox{\boldmath $R$}\,dx
\\
\\
\displaystyle
&
\displaystyle
\,\,\,
+\rho \theta_0\tau e^{-\tau T}\int_{\Omega\setminus\overline D}(\mbox{\boldmath $F$}
-
\mbox{\boldmath $F$}_0)\cdot\mbox{\boldmath $R$}\,dx.
\end{array}
\tag {2.21}
$$
Summing both sides of (2.20) and (2.21), we obtain
$$\begin{array}{l}
\displaystyle
\,\,\,\,\,\,
\int_{\partial D}k\frac{\partial\Xi_0}{\partial\mbox{\boldmath $\nu$}}\Sigma\,dS
+\theta_0\tau\int_{\partial D}s(\mbox{\boldmath $w$}_0,\Xi_0)\mbox{\boldmath $\nu$}\cdot
\mbox{\boldmath $R$}\,dS
\\
\\
\displaystyle
=e(\tau)+\theta_0\tau E(\tau)\\
\\
\displaystyle
\,\,\,
+\rho e^{-\tau T}\int_{\Omega\setminus\overline D}(h-h_0)\Sigma\,dx
+\rho \theta_0\tau e^{-\tau T}\int_{\Omega\setminus\overline D}(\mbox{\boldmath $F$}
-
\mbox{\boldmath $F$}_0)\cdot\mbox{\boldmath $R$}\,dx.
\end{array}
\tag{2.22}
$$
Rewrite this right-hand side as
$$\begin{array}{l}
\displaystyle
\,\,\,\,\,\,
\int_{\Omega\setminus\overline D}
\left(k\vert\nabla\Sigma\vert^2
+c\tau\left\vert
\Sigma+\frac{\rho e^{-\tau T}}{2c\tau}(h-h_0)\right\vert^2\right)\,dx
\\
\\
\displaystyle
\,\,\,
+\theta_0\tau
\int_{\Omega\setminus\overline D}
\left(2\mu\left\vert\text{Sym}\,\nabla\mbox{\boldmath $R$}\right\vert^2+\lambda\vert\nabla\cdot\mbox{\boldmath $R$}\vert^2
+\rho\tau^2
\left\vert\mbox{\boldmath $R$}+\frac{e^{-\tau T}}{2\tau^2}(\mbox{\boldmath $F$}-\mbox{\boldmath $F$}_0)\right\vert^2\right)\,dx\\
\\
\displaystyle
-\frac{\rho^2 e^{-2\tau T}}{4c\tau}\Vert h-h_0\Vert^2_{L^2(\Omega\setminus\overline D)}
-\frac{\rho\theta_0 e^{-2\tau T}}{4\tau}\Vert\mbox{\boldmath $F$}-\mbox{\boldmath $F$}_0\Vert^2_{L^2(\Omega\setminus\overline D)}.
\end{array}
$$
Since $\Vert h-h_0\Vert_{L^2(\Omega)}=O(1)$ and $\Vert\mbox{\boldmath $F$}-\mbox{\boldmath $F$}_0\Vert_{L^2(\Omega\setminus\overline D)}=O(\tau)$,
it follows from this expression and (2.22) that
$$
\displaystyle
e(\tau)+\theta_0\tau E(\tau)
\le
\int_{\partial D}k\frac{\partial\Xi_0}{\partial\mbox{\boldmath $\nu$}}\Sigma\,dS
+\theta_0\tau\int_{\partial D}s(\mbox{\boldmath $w$}_0,\Xi_0)\mbox{\boldmath $\nu$}\cdot
\mbox{\boldmath $R$}\,dS+O(\tau e^{-2\tau T}).
\tag {2.23}
$$
Take a lifting  $\tilde{\Sigma}$ of $\Sigma\vert_{\partial D}$ into $D$ such that $\Vert\tilde{\Sigma}\Vert_{H^1(D)}
\le C\Vert\Sigma\Vert_{H^1(\Omega\setminus\overline D)}$, where $C$ is a positive constant independent of $\Sigma$.
Then, the second equation on (2.3) gives
$$\begin{array}{ll}
\displaystyle
\int_{\partial D}k\frac{\partial\Xi_0}{\partial\mbox{\boldmath $\nu$}}\Sigma\,dS
&
\displaystyle
=\int_{D}
(c\tau\Xi_0-m\theta_0\tau\nabla\cdot\mbox{\boldmath $w$}_0)\tilde{\Sigma}\,dx
+\int_D k\nabla\Xi_0\cdot\nabla\tilde{\Sigma}\,dx\\
\\
\displaystyle
\,\,\,
&
\displaystyle
\,\,\,
+e^{-\tau T}\int_D h_0\tilde{\Sigma}\,dx.
\end{array}
\tag {2.24}
$$
Here from (2.8), (2.14), (2.15) we have
$$\left\{
\begin{array}{l}
\Vert\Sigma\Vert _{H^1(\Omega\setminus\overline D)}=O(e(\tau)^{1/2}),\\
\\
\displaystyle
\Vert\Xi_0\Vert_{L^2(D)}=O(\tau^{-1/2}j(\tau)^{1/2}),
\\
\\
\displaystyle
\Vert\nabla\Xi_0\Vert_{L^2(D)}=O(j(\tau)^{1/2}),
\\
\\
\displaystyle
\Vert\nabla\cdot\mbox{\boldmath $w$}_0\Vert_{L^2(D)}
=O(J(\tau)^{1/2}).
\end{array}
\right.
$$
Applying these to the right-hand side on (2.24),
we obtain
$$\begin{array}{l}
\,\,\,\,\,\,
\displaystyle
\int_{\partial D}k\frac{\partial\Xi_0}{\partial\mbox{\boldmath $\nu$}}\Sigma\,dS
\le C e(\tau)^{1/2}(\tau^{1/2}j(\tau)^{1/2}+\tau J(\tau)^{1/2}+e^{-\tau T}).
\end{array}
\tag {2.25}
$$
A similar technique together with the first equation on (2.3) gives also
$$\begin{array}{l}
\displaystyle
\,\,\,\,\,\,
\theta_0\tau\int_{\partial D}s(\mbox{\boldmath $w$}_0,\Xi_0)\mbox{\boldmath $\nu$}\cdot
\mbox{\boldmath $R$}\,dS\\
\\
\displaystyle
\le C\tau\Vert\mbox{\boldmath $R$}\Vert_{H^1(\Omega\overline D)}
\left(\tau^2\Vert\mbox{\boldmath $w$}_0\Vert_{L^2(D)}
+\Vert\nabla\mbox{\boldmath $w$}_0\Vert_{L^2(D)}+
\Vert\nabla\Xi_0\Vert_{L^2(D)}
+\tau e^{-\tau T}\right).
\end{array}
$$
Korn's second inequality \cite{DuL} tells us that
$$
\left\{
\begin{array}{l}
\displaystyle
\Vert\mbox{\boldmath $w$}_0\Vert_{H^1(D)}
\le C''
\left(\Vert\mbox{\boldmath $w$}_0\Vert_{L^2(D)}^2
+\Vert\text{Sym}\,\nabla\mbox{\boldmath $w$}_0\Vert_{L^2(D)}^2\right)^{1/2},
\\
\\
\displaystyle
\Vert\mbox{\boldmath $R$}\Vert_{H^1(\Omega\setminus\overline D)}
\le C''
\left(\Vert\mbox{\boldmath $R$}\Vert_{L^2(\Omega\setminus\overline D)}^2
+\Vert\text{Sym}\,\nabla\mbox{\boldmath $R$}\Vert_{L^2(\Omega\setminus\overline D)}^2\right)^{1/2},
\end{array}
\right.
\tag {2.26}
$$
where $C''$ is a positive constant independent of $\mbox{\boldmath $w$}_0$ and $\mbox{\boldmath $R$}$.
Here we note that, for all real $3\times 3$-matrix $A$ we have
$$\begin{array}{l}
\displaystyle
\,\,\,\,\,\,
\vert\text{Sym}\,A\vert^2
\\
\\
\displaystyle
\le
2\left(\left\vert\text{Sym}\,A-\frac{\text{trace}\,A}{3}I_3\right\vert^2+\frac{\vert\text{trace}\,A\vert^2}{3}
\right)
\\
\\
\displaystyle
=\frac{2}{2\mu}\cdot 2\mu\left\vert\text{Sym}\,A-\frac{\text{trace}\,A}{3}I_3\right\vert^2+
\frac{2}{3\lambda+2\mu}\cdot\frac{3\lambda+2\mu}{3}\vert\text{trace}\,A\vert^2\\
\\
\displaystyle
\le 2\max\left\{\frac{1}{2\mu}, \frac{1}{3\lambda+2\mu}\right\}
\left(
2\mu\left\vert\text{Sym}\,A-\frac{\text{trace}\,A}{3}I_3\right\vert^2+
\frac{3\lambda+2\mu}{3}\vert\text{trace}\,A\vert^2
\right).
\end{array}
$$
Applying identity (2.12) to this right-hand side,
we obtain
$$\displaystyle
\vert\text{Sym}\,A\vert^2
\le C
\left(2\mu\vert\text{Sym}\,A\vert^2+\lambda\vert\text{trace}\,A\vert^2\right)
\tag {2.27}
$$
provided $\mu>0$ and $3\lambda+2\mu>0$ with a suitable positive constant $C=C(2\mu,3\lambda+2\mu)$.

Thus, substituting $A=\nabla\mbox{\boldmath $w$}_0, \nabla\mbox{\boldmath $R$}$ into (2.27) and using (2.26),
repectively,
for all $\tau$ with $\rho\tau^2\ge 1$ we have
$$
\left\{
\begin{array}{l}
\displaystyle
\Vert\mbox{\boldmath $w$}_0\Vert_{H^1(D)}
=O(J(\tau)^{1/2}),\\
\\
\displaystyle
\Vert\mbox{\boldmath $R$}\Vert_{H^1(\Omega\setminus\overline D)}
=O(E(\tau)^{1/2}).
\end{array}
\right.
\tag {2.28}
$$
We have also the following trivial estimates:
$$
\left\{
\begin{array}{l}
\displaystyle
\Vert\mbox{\boldmath $w$}_0\Vert_{L^2(D)}=O(\tau^{-1}J(\tau)^{1/2}),\\
\\
\displaystyle
\Vert\nabla\Xi_0\Vert_{L^2(D)}=O(j(\tau)^{1/2}).
\end{array}
\right.
$$
Hence we obtain
$$\displaystyle
\theta_0\tau\int_{\partial D}s(\mbox{\boldmath $w$}_0,\Xi_0)\mbox{\boldmath $\nu$}\cdot
\mbox{\boldmath $R$}\,dS
\le C^{'''}\tau E(\tau)^{1/2}(\tau J(\tau)^{1/2}+j(\tau)^{1/2}+\tau e^{-\tau T}).
\tag {2.29}
$$
Summing both sides on (2.25) and (2.29) and then from (2.23) we obtain
$$\begin{array}{ll}
\displaystyle
\,\,\,\,\,\,
e(\tau)+\theta_0\tau E(\tau)
&
\displaystyle
\le C_1 e(\tau)^{1/2}(\tau^{1/2}j(\tau)^{1/2}+\tau J(\tau)^{1/2}+e^{-\tau T})\\
\\
\displaystyle
&
\,\,\,
+C_2\tau E(\tau)^{1/2}(\tau J(\tau)^{1/2}+j(\tau)^{1/2}+\tau e^{-\tau T})
+C_3\tau e^{-2\tau T}.
\end{array}
$$
Now a standard technique gives (2.19).

\noindent
$\Box$

\subsection{Three special solutions and their bounds}

In this subsection we introduce three functions and describe their upper and lower bounds.

Let $\mbox{\boldmath $w$}_{s0}\in H^1(\Bbb R^3)^3$ solve
$$\begin{array}{ll}
\displaystyle
\mu\Delta\mbox{\boldmath $w$}_{s0}-\rho\tau^2\mbox{\boldmath $w$}_{s0}
+\rho\mbox{\boldmath $v$}_0(x)=\mbox{\boldmath $0$} & \text{in $\Bbb R^3$,}
\end{array}
$$
where $\mbox{\boldmath $v$}_0(x)=\partial_t\mbox{\boldmath $v$}_s(x,0)$ and $\mbox{\boldmath $v$}_s$ is giveny by (1.10).
Since we have
$$\displaystyle
\mbox{\boldmath $v$}_0(x)=-2\chi_B(x)(x-p)\times\mbox{\boldmath $a$},
$$
$\mbox{\boldmath $w$}_{s0}$ takes the form
$$\begin{array}{ll}
\displaystyle
\mbox{\boldmath $w$}_{s0}(x)
&
\displaystyle
=-\frac{\rho}{2\pi}
\left(
\int_B\frac{e^{-\tau\sqrt{\rho/\mu}\,\vert x-y\vert}}{\vert x-y\vert}\,(y-p)dy\right)\times\mbox{\boldmath $a$}
\end{array}
\tag {2.30}
$$
and we have $\nabla\cdot\mbox{\boldmath $w$}_{s0}=0$.

Second, let $\mbox{\boldmath $w$}_{p0}\in H^1(\Bbb R^3)^3$ solve
$$\begin{array}{ll}
\displaystyle
(\lambda+2\mu)\Delta\mbox{\boldmath $w$}_{p0}-\rho\tau^2\mbox{\boldmath $w$}_{p0}
+\rho\mbox{\boldmath $v$}_0(x)=\mbox{\boldmath $0$} & \text{in $\Bbb R^3$,}
\end{array}
$$
where $\mbox{\boldmath $v$}_0(x)=\partial_t\mbox{\boldmath $v$}_p(x,0)$ and $\mbox{\boldmath $v$}_p$
is given by (1.15).
Since we have
$$\displaystyle
\mbox{\boldmath $v$}_0(x)=-2\chi_B(x)(\eta-\vert x-p\vert)\frac{x-p}{\vert x-p\vert},
$$
one gets the expression
$$\begin{array}{ll}
\displaystyle
\mbox{\boldmath $w$}_{p0}(x)
&
\displaystyle
=-\frac{\rho}{2\pi}
\int_B\frac{e^{-\tau\sqrt{\rho/(\lambda+2\mu)}\,\vert x-y\vert}}{\vert x-y\vert}\,(\eta-\vert y-p\vert)\frac{y-p}{\vert y-p\vert}\,dy
\end{array}
\tag {2.31}
$$
and we have $\nabla\times\mbox{\boldmath $w$}_{p0}=\mbox{\boldmath $0$}$.

Third, let $\Theta_{00}\in H^1(\Bbb R^3)$ solve
$$\begin{array}{ll}
\displaystyle
k\Delta\Theta_{00}-c\tau\Theta_{00}+cf_0=0 & \text{in $\Bbb R^3$,}
\end{array}
$$
where $f_0(x)=\Theta_0(x,0)$ and $\Theta_0$ is the solution of (1.16).
Then, we have the expression
$$\begin{array}{ll}
\displaystyle
\Theta_{00}(x)
&
\displaystyle
=\frac{c}{4\pi k}
\int_B\frac{e^{-\sqrt{\tau}\sqrt{c/k}\,\vert x-y\vert}}{\vert x-y\vert}\,(\eta-\vert y-p\vert)^2\,dy.
\end{array}
\tag {2.32}
$$
It is easy to see that from the expression (2.30) to (2.32) one gets the following estimates.

\proclaim{\noindent Lemma 2.1.}
Let $U$ be an arbitrary bounded open subset of $\Bbb R^3$ such that $\overline B\cap\overline U=\emptyset$.
Then, we have, as $\tau\longrightarrow\infty$
$$
\left\{
\begin{array}{l}
\displaystyle
\tau
\Vert
\mbox{\boldmath $w$}_{s0}
\Vert_{L^2(U)}
+
\Vert\nabla
\mbox{\boldmath $w$}_{s0}
\Vert_{L^2(U)}
=O(\tau e^{-\tau\sqrt{\rho/\mu}\,\text{dist}\,(U,B)}),\\
\\
\displaystyle
\tau
\Vert
\mbox{\boldmath $w$}_{p0}
\Vert_{L^2(U)}
+
\Vert\nabla
\mbox{\boldmath $w$}_{p0}
\Vert_{L^2(U)}
=O(\tau e^{-\tau\sqrt{\rho/(\lambda+2\mu)}\,\text{dist}\,(U,B)}),
\\
\\
\displaystyle
\sqrt{\tau}
\Vert
\Theta_{00}
\Vert_{L^2(U)}
+
\Vert\nabla
\Theta_{00}
\Vert_{L^2(U)}
=O(\sqrt{\tau}\,e^{-\sqrt{\tau}\,\sqrt{c/k}\,\text{dist}\,(U,B)}).
\end{array}
\right.
\tag {2.33}
$$

\endproclaim

The estimates in the next two lemmas play the key role in this paper.

\proclaim{\noindent Lemma 2.2.}
Let $R>\eta$.  Then, there exist positive numbers $C$ and $\tau_0$ such that, for all $\tau\ge\tau_0$
and $x\in\Bbb R^3\setminus\overline D$ with $\vert x-p\vert\le R$
$$\displaystyle
\vert\mbox{\boldmath $w$}_{so}(x)\vert
\ge C\tau^{-1}e^{-\tau\sqrt{\rho/\mu}\,(\vert x-p\vert-\eta)}\,
\left\vert\frac{x-p}{\vert x-p\vert}\times\mbox{\boldmath $a$}\right\vert
\tag {2.34}
$$
and
$$\displaystyle
\vert\mbox{\boldmath $w$}_{po}(x)\vert
\ge C\tau^{-1}e^{-\tau\sqrt{\rho/(\lambda+2\mu)}\,(\vert x-p\vert-\eta)}.
\tag {2.35}
$$

\endproclaim
{\it\noindent Proof.}
From (2.30) and (2.31) one gets
$$\begin{array}{ll}
\displaystyle
\mbox{\boldmath $w$}_{s0}(x)
&
\displaystyle
=-\frac{\rho}{2\pi}\mbox{\boldmath $I$}_1(x;\tau\sqrt{\rho/\mu})\times\mbox{\boldmath $a$}
\end{array}
\tag {2.36}
$$
and
$$\begin{array}{ll}
\displaystyle
\mbox{\boldmath $w$}_{p0}(x)
&
\displaystyle
=-\frac{\rho}{2\pi}
(\eta\mbox{\boldmath $I$}_0(x;\tau\sqrt{\rho/(\lambda+2\mu)})-\mbox{\boldmath $I$}_1(x;\tau\sqrt{\rho/(\lambda+2\mu)})),
\end{array}
\tag {2.37}
$$
where
$$\left\{
\begin{array}{l}
\displaystyle
\mbox{\boldmath $I$}_0(x;\tau)=\int_B\frac{e^{-\tau\,\vert x-y\vert}}{\vert x-y\vert}\frac{y-p}{\vert y-p\vert}\,dy,
\\
\\
\displaystyle
\mbox{\boldmath $I$}_1(x;\tau)=\int_B\frac{e^{-\tau\,\vert x-y\vert}}{\vert x-y\vert}(y-p)\,dy.
\end{array}
\right.
\tag {2.38}
$$
It follows from the middle equation on (A.3) and (A.5) that
$$\begin{array}{l}
\displaystyle
\mbox{\boldmath $I$}_1(x;\tau)
\sim
\frac{4\pi}{\tau^2}
\left(\vert x-p\vert^2\cdot\tau\eta\cdot\frac{e^{\tau\eta}}{2}
+\frac{1}{\tau^2}\cdot(\tau\eta)^3\cdot\frac{e^{\tau\eta}}{2}
\right)
\cdot\frac{e^{-\tau\vert x-p\vert}}{\vert x-p\vert^2}\cdot
\frac{x-p}{\vert x-p\vert}\\
\\
\displaystyle
=\frac{2\pi\eta}{\tau}(\vert x-p\vert^2+\eta^2)e^{\tau\eta}\frac{x-p}{\vert x-p\vert}.
\end{array}
$$
A combination of this and (2.36) yields (2.34).

(2.35) is proved as follows.
From the first equation on (A.3) and (A.4) we have
$$\begin{array}{l}
\displaystyle
\eta\mbox{\boldmath $I$}_0(x;\tau)
\sim
\pi\left(\frac{2}{\tau}\cdot\vert x-p\vert^2\cdot\frac{e^{\tau\eta}}{2}
+\frac{2}{\tau^3}\cdot(\tau\eta)^2\cdot\frac{e^{\tau\eta}}{2}
\right)
\frac{e^{-\tau\vert x-p\vert}}{\vert x-p\vert^2}\cdot\frac{x-p}{\vert x-p\vert}
\\
\\
\displaystyle
=\frac{\pi}{\tau}(\vert x-p\vert^2+\eta^2)\cdot\frac{e^{-\tau\,(\vert x-p\vert-\eta)}}{\vert x-p\vert^2}\cdot\frac{x-p}{\vert x-p\vert}.
\end{array}
$$
From the second equation on (A.3) and (A.5) we have
$$\begin{array}{l}
\displaystyle
\mbox{\boldmath $I$}_1(x;\tau)
\sim
\frac{4\pi}{\tau^2}
\left(\vert x-p\vert^2\cdot\tau\eta\cdot\frac{e^{\tau\eta}}{2}
+\frac{1}{\tau^2}\cdot(\tau\eta)^3\cdot\frac{e^{\tau\eta}}{2}
\right)
\frac{e^{-\tau\vert x-p\vert}}{\vert x-p\vert^2}\cdot\frac{x-p}{\vert x-p\vert}
\\
\\
\displaystyle
=\frac{2\pi}{\tau}(\vert x-p\vert^2+\eta^2)\cdot\frac{e^{-\tau\,(\vert x-p\vert-\eta)}}{\vert x-p\vert^2}\cdot\frac{x-p}{\vert x-p\vert}.
\end{array}
$$
Thus, we obtain
$$\displaystyle
\eta\mbox{\boldmath $I$}_0(x;\tau)
-\mbox{\boldmath $I$}_1(x;\tau)
\sim
-\frac{\pi}{\tau}(\vert x-p\vert^2+\eta^2)\cdot\frac{e^{-\tau\,(\vert x-p\vert-\eta)}}{\vert x-p\vert^2}\cdot\frac{x-p}{\vert x-p\vert}.
$$
Then (2.37) yields the desired estimate.

\noindent
$\Box$

\proclaim{\noindent Lemma 2.3.}
There exist positive constants $C$ and $\tau_0$ such that fo all $\tau\ge\tau_0$ and all $x\in\Bbb R^3\setminus\overline B$,
$$\displaystyle
\Theta_{00}(x)\ge C\tau^{-3/2}\frac{e^{-\sqrt{\tau}\sqrt{c/k}\,(\vert x-p\vert-\eta)}}{\vert x-p\vert}.
\tag {2.39}
$$
\endproclaim
{\it\noindent Proof.}
From (2.32) one has the expression
$$\begin{array}{ll}
\displaystyle
\Theta_{00}(x)
&
\displaystyle
=\frac{c}{4\pi k}
(\eta^2 I_0(x;\sqrt{\tau}\sqrt{c/k})-2\eta I_1(x;\sqrt{\tau}\sqrt{c/k})+I_2(x;\sqrt{\tau}\sqrt{c/k})),
\end{array}
\tag {2.40}
$$
where
$$\displaystyle
I_j(x;\tau)
=\int_B\frac{e^{-\tau\,\vert x-y\vert}}{\vert x-y\vert}\,\vert y-p\vert^j\,dy.
\tag {2.41}
$$
From (A.1), (A.2), the last equation on (A.3) and (A.6) we have:
$$\displaystyle
\eta^2I_0(x;\tau)
\sim
\frac{4\pi\eta^2}{\tau^3}\cdot
(\tau\eta-1)\frac{e^{\tau\eta}}{2}
\cdot\frac{e^{-\tau\vert x-p\vert}}{\vert x-p\vert}
=\frac{2\pi\eta^2}{\tau^3}(\tau\eta-1)
\cdot
\frac{e^{-\tau(\vert x-p\vert-\eta)}}{\vert x-p\vert};
$$
$$
\displaystyle
2\eta I_1(x;\tau)
\sim
\frac{8\pi\eta}{\tau^4}
\cdot
\tau\eta(\tau\eta-2)
\frac{e^{\tau\eta}}{2}
\cdot\frac{e^{-\tau\vert x-p\vert}}{\vert x-p\vert}
=\frac{4\pi\eta^2}{\tau^3}(\tau\eta-2)\cdot\frac{e^{-\tau(\vert x-p\vert-\eta)}}{\vert x-p\vert};
$$
$$\displaystyle
I_2(x;\tau)
\sim\frac{4\pi}{\tau^5}\cdot
(\tau\eta)^2(\tau\eta+6-3)
\frac{e^{\tau\eta}}{2}
\cdot
\frac{e^{-\tau\vert x-p\vert}}{\vert x-p\vert}
=\frac{2\pi\eta^2}{\tau^3}\cdot
(\tau\eta+3)
\frac{e^{-\tau(\vert x-p\vert-\eta)}}{\vert x-p\vert}.
$$
Thus we have
$$\displaystyle
\eta^2I_0(x;\tau)
-2\eta I_1(x;\tau)
+I_2(x;\tau)
\sim \frac{12\pi\eta^2}{\tau^3}\frac{e^{-\tau(\vert x-p\vert-\eta)}}{\vert x-p\vert}.
$$
A combination of this and (2.40) gives the desired conclusion.

\noindent
$\Box$

\subsection{Bounds for $J(\tau)$ and $j(\tau)$ for special $\mbox{\boldmath $v$}$ and $\Theta$}

The integrals $J(\tau)$ and $j(\tau)$ given by (2.8) and (2.14), respectively,
depends on $m$ and the pair $\mbox{\boldmath $v$}$ and $\Theta$ which is a solution
of (1.2) with (1.3).
To make it clear, in this subsection, we denote $J(\tau)$ and $j(\tau)$ by $J_m(\tau;\mbox{\boldmath $v$},\Theta)$
and $j_m(\tau;\mbox{\boldmath $v$},\Theta)$, repectively.
The same remark works also for the pair of $\mbox{\boldmath $w$}_0$ and $\Xi_0$ which are given by the second quations on (1.8) and (1.9).
We denote them by $\mbox{\boldmath $w$}_0^m(\,\cdot\,,\mbox{\boldmath $v$},\Theta)$ and 
$\Xi_0^m(\,\cdot\,,\mbox{\boldmath $v$},\Theta)$, respectively.

\proclaim{\noindent Proposition 2.4.}
Let $U$ be a bounded open subset of $\Bbb R^3$ with $\overline B\cap\overline U=\emptyset$.
Then, we have, as $\tau\longrightarrow\infty$
$$
\displaystyle
\tau\Vert\mbox{\boldmath $w$}_0^m(\,\cdot\,,\mbox{\boldmath $v$}_s,0)\Vert_{L^2(U)}+
\Vert\nabla\mbox{\boldmath $w$}_0^m(\,\cdot\,,\mbox{\boldmath $v$}_s,0)\Vert_{L^2(U)}
=O(\tau e^{-\tau\sqrt{\rho/\mu}\,\text{dist}\,(U,B)}+e^{-\tau T}),
\tag {2.42}
$$
$$
\displaystyle
\tau\Vert\mbox{\boldmath $w$}_0^0(\,\cdot\,,\mbox{\boldmath $v$}_p,0\Vert_{L^2(U)}+
\Vert\nabla\mbox{\boldmath $w$}_0^0(\,\cdot\,,\mbox{\boldmath $v$}_p,0)\Vert_{L^2(U)}
=O(\tau e^{-\tau\sqrt{\rho/(\lambda+2\mu)}\,\text{dist}\,(U,B)}+e^{-\tau T}),
\tag {2.43}
$$
$$
\displaystyle
\sqrt{\tau}\Vert\Xi_0^0(\,\cdot\,,\mbox{\boldmath $0$},\Theta_0)\Vert_{L^2(U)}+
\Vert\nabla\Xi_0^0(\,\cdot\,,\mbox{\boldmath $0$},\Theta_0)\Vert_{L^2(U)}
=O(\sqrt{\tau}e^{-\sqrt{\tau}\sqrt{c/k}\,\text{dist}\,(U,B)}+e^{-\tau T}).
\tag {2.44}
$$

\endproclaim

{\it\noindent Proof.}
First we give a proof of (2.42).
Set
$$
\displaystyle
\mbox{\boldmath $\epsilon$}_s=e^{\tau T}(\mbox{\boldmath $w$}_0-\mbox{\boldmath $w$}_{s0}).
$$
We have 
$$\displaystyle
\mbox{\boldmath $w$}_0=\mbox{\boldmath $w$}_{s0}+e^{-\tau T}
\mbox{\boldmath $\epsilon$}_s
\tag {2.45}
$$
and $\nabla\cdot\mbox{\boldmath $w$}_0=0$ since $\nabla\cdot\mbox{\boldmath $v$}_s=0$.
These together with $\nabla\cdot\mbox{\boldmath $w$}_{s0}=0$ yield
$\nabla\cdot\mbox{\boldmath $\epsilon$}_s=0$.
Moreover, we have $\Xi_0=0$.
Then, from the first equation on (2.3) we have
$$\begin{array}{ll}
\displaystyle
(\mu\Delta-\rho\tau^2)\mbox{\boldmath $\epsilon$}_s=\rho\mbox{\boldmath $F$}_0 & \text{in}\,\Bbb R^3.
\end{array}
\tag {2.46}
$$
Then, using the third equation on (2.4), we can easily see that 
$$
\displaystyle \tau\Vert\mbox{\boldmath $\epsilon$}_s\Vert_{L^2(\Bbb R^3)}+
\Vert\nabla\mbox{\boldmath $\epsilon$}_s\Vert_{L^2(\Bbb R^3)}=O(1).
\tag {2.47}
$$
A combination of this and the first estimate on (2.33) yields (2.42).

Next set
$$
\displaystyle
\mbox{\boldmath $\epsilon$}_p=e^{\tau T}(\mbox{\boldmath $w$}_0-\mbox{\boldmath $w$}_{p0}).
$$
We have
$$\displaystyle
\mbox{\boldmath $w$}_0=\mbox{\boldmath $w$}_{p0}+e^{-\tau T}
\mbox{\boldmath $\epsilon$}_p
$$
and $\nabla\times\mbox{\boldmath $w$}_0=\mbox{\boldmath $0$}$ since $\nabla\times\mbox{\boldmath $v$}_p=\mbox{\boldmath $0$}$.
These together with $\nabla\times\mbox{\boldmath $w$}_{p0}=\mbox{\boldmath $0$}$ yield
$\nabla\times\mbox{\boldmath $\epsilon$}_p=\mbox{\boldmath $0$}$.
Then, from the first equation on (2.3) with $m=0$ and the equation 
$\nabla(\nabla\cdot\mbox{\boldmath $\epsilon$}_p)=\Delta\mbox{\boldmath $\epsilon$}_p$,
we have
$$\begin{array}{ll}
\displaystyle
\{(\lambda+2\mu)\Delta-\rho\tau^2\}\mbox{\boldmath $\epsilon$}_p=\rho\mbox{\boldmath $F$}_0 & \text{in}\,\Bbb R^3.
\end{array}
$$
Then, we have 
$$
\displaystyle \tau\Vert\mbox{\boldmath $\epsilon$}_p\Vert_{L^2(\Bbb R^3)}+
\Vert\nabla\mbox{\boldmath $\epsilon$}_p\Vert_{L^2(\Bbb R^3)}=O(1).
$$
A combination of this and the middle estimate on (2.33) yields (2.43).
(2.44) is also a consequence of the last estimate on (2.33) and the second equation on (2.3) with $m=0$.

\noindent
$\Box$

\proclaim{\noindent Proposition 2.5.}
(i)  We have
$$\left\{
\begin{array}{ll}
\displaystyle
J_m(\tau;\mbox{\boldmath $v$}_s,0)
\displaystyle
=O(\tau^2 e^{-2\tau\sqrt{\rho/\mu}\,\text{dist}\,(D,B)}+e^{-2\tau T})
&
\text{as $\tau\longrightarrow\infty$,}
\\
\\
\displaystyle
j_m(\tau;\mbox{\boldmath $v$}_s,0)=0
&
\text{for all $\tau>0$.}
\end{array}
\right.
\tag {2.48}
$$
(ii)  We have, as $\tau\longrightarrow\infty$
$$\left\{\begin{array}{l}
\displaystyle
J_0(\tau;\mbox{\boldmath $v$}_p,0)
\displaystyle
=O(\tau^2 e^{-2\tau\sqrt{\rho/(\lambda+2\mu)}\,\text{dist}\,(D,B)}+e^{-2\tau T}),
\\
\\
\displaystyle
j_0(\tau;\mbox{\boldmath $0$},\Theta_0)=O(\tau e^{-2\sqrt{\tau}\,\sqrt{c/k}\,\text{dist}\,(D,B)}+e^{-2\tau T}).
\end{array}
\right.
$$
(iii) Let $T$ satisfies 
$$\displaystyle
T>\sqrt{\frac{\rho}{\mu}}\,\text{dist}\,(D,B).
\tag {2.49}
$$
Then, there exist positive constants $\tau_0$ and $C$ such that, for all $\tau\ge\tau_0$
$$\displaystyle
\tau^{3} e^{2\tau\sqrt{\rho/\mu}\,\text{dist}\,(D,B)}J_m(\tau;\mbox{\boldmath $v$}_s,0)\ge C.
\tag {2.50}
$$
(iv) Let $T$ satisfies 
$$\displaystyle
T>\sqrt{\frac{\rho}{\lambda+2\mu}}\,\text{dist}\,(D,B).
$$
Then, there exist positive constants $\tau_0$ and $C$ such that, for all $\tau\ge\tau_0$
$$\displaystyle
\tau^{2} e^{2\tau\sqrt{\rho/(\lambda+2\mu)}\,\text{dist}\,(D,B)}J_0(\tau;\mbox{\boldmath $v$}_p,0)\ge C.
$$
(v) Let $T$ be an arbitrary positive number.
Then, there exist positive constants $\tau_0$ and $C$ such that, for all $\tau\ge\tau_0$
$$\displaystyle
\tau^{4} e^{2\sqrt{\tau}\,\sqrt{c/k}\,\text{dist}\,(D,B)}j_0(\tau;\mbox{\boldmath $0$},\Theta_0)\ge C.
$$

\endproclaim

{\it\noindent Proof.}
Applying Proposition 2.4 in the case when $U=D$ to the expression (2.8) and (2.14) for $J(\tau)$ and $j(\tau)$,
we obtain the first estimate on (2.48) and (ii).  From $\Xi_0^m(\,\cdot\,;\mbox{\boldmath $v$}_s,0)=0$, we
obtain $j_m(\tau;\mbox{\boldmath $v$}_s,0)=0$.

Next we give the proof of (iii).
It follows from (2.8) and (2.45) and (2.47)
$$
\displaystyle
J_m(\tau;\mbox{\boldmath $v$}_s,0)
\ge \frac{1}{2}J_0(\tau)+O(e^{-2\tau T}),
\tag {2.51}
$$
where
$$\displaystyle
J_0(\tau)=\int_D(2\mu\left\vert\text{Sym}\,\nabla\mbox{\boldmath $w$}_{s0}\right\vert^2+\rho\tau^2\vert\mbox{\boldmath $w$}_{s0}\vert^2)dx.
$$
From (2.34) in Lemma 2.2 we have
$$\displaystyle
J_0(\tau)
\ge
C^2\int_D
e^{-2\tau\sqrt{\rho/\mu}\,(\vert x-p\vert-\eta)}
\left\vert\frac{x-p}{\vert x-p\vert}\times\mbox{\boldmath $a$}\right\vert^2
\,dx.
$$
Applying Lemma A.1, we conclude that
here exist positive constants $\tau_0$, $C'$  and $\kappa$ such that,
for all $\tau\ge\tau_0$ 
$$\displaystyle
\tau^{\kappa}e^{2\tau\sqrt{\rho/\mu}\,\text{dist}\,(D,B)}\,J_0(\tau)\ge C'.
$$
Now from this and (2.51)  we see that  (2.50) is valid under condition (2.49).
Note that $\kappa=2$ or $\kappa=3$.  Thus we have chosen the worst $\kappa$.

(iv) and (v) are easy consequences of (2.35) in Lemma 2.2, (2.39) in Lemma 2.3, Lemma A.2.

\noindent
$\Box$

{\bf\noindent Remark 2.2.}
(ii), (iv) and (v) are for the proof of Theorems 1.2 and 1.3.

\section{Proof of Theorem 1.1.}

Once we have the following lower and upper bounds for indicator function
$I^1(\tau;\mbox{\boldmath $v$}_s,0)$ as $\tau\longrightarrow\infty$, then the proof of Theorem 1.1
is a due cource as we have done in \cite{IE4}.  Just apply the first estimate on (2.48) and (iii) in Proposition 2.5.

\proclaim{\noindent Proposition 3.1.}
Let $T$ be an arbitraly positive number.
We have, as $\tau\longrightarrow\infty$
$$\displaystyle
I^1(\tau;\mbox{\boldmath $v$}_s,0)
=O(\tau^2 J_m(\tau;\mbox{\boldmath $v$}_s,0)
+\tau^2e^{-\tau T}e^{-\tau\sqrt{\rho/\mu}\,\text{dist}\,(\Omega,B)}
+\tau^2 e^{-2\tau T})
\tag {3.1}
$$
and
$$\displaystyle
I^1(\tau;\mbox{\boldmath $v$}_s,0)
\ge J_m(\tau;\mbox{\boldmath $v$}_s,0)
+O(\tau^2e^{-\tau T}e^{-\tau\sqrt{\rho/\mu}\,\text{dist}\,(\Omega,B)}
+\tau^2 e^{-2\tau T}).
\tag {3.2}
$$
\endproclaim

{\it\noindent Proof.}
As pointed out in Section 1.1, we have
$\Xi_0=0$ and hence
$$
\displaystyle
I^2(\tau;\mbox{\boldmath $v$}_s,0)=0.
$$
Note also that a combination of the second equation on (1.8) and (1.10) gives
$$\displaystyle
\nabla\cdot\mbox{\boldmath $w$}_0=0.
$$
Thus from (2.17) we obtain
$$
\displaystyle
I^1(\tau;\mbox{\boldmath $v$}_s,0)
=J_m(\tau;\mbox{\boldmath $v$}_s,0)+
\left(E(\tau)
+\frac{1}{\theta_0\tau}
e(\tau)\right)
+{\cal R}(\tau),
\tag {3.3}
$$
where ${\cal R}(\tau)$ is given by (2.18).

A combination of (2.19) and the second equation on (2.48)
gives
$$
\displaystyle
E(\tau)+\frac{1}{\theta_0\tau}e(\tau)
=O(\tau^2J_m(\tau;\mbox{\boldmath $v$}_s,0)+\tau^2 e^{-2\tau T}).
\tag {3.4}
$$
Thus, for the proof of (3.1) and (3.2) it sfficies to give an estimate of ${\cal R}(\tau)$.
We show that, as $\tau\longrightarrow\infty$
$$\displaystyle
{\cal R}(\tau)
=O(\tau^2e^{-\tau T}e^{-\tau\sqrt{\rho/\mu}\,\text{dist}\,(\Omega,B)}
+\tau e^{-2\tau T}).
\tag {3.5}
$$
A combination of (3.4) and the first estimate on (2.48) gives
$$
\displaystyle
E(\tau)+\frac{1}{\theta_0\tau}e(\tau)
=O(\tau^4 e^{-2\tau\sqrt{\rho/\mu}\,\text{dist}\,(D,B)}+\tau^2 e^{-2\tau T}).
\tag {3.6}
$$
This and (2.9) yields
$$\displaystyle
\Vert\mbox{\boldmath $R$}\Vert_{L^2(\Omega\setminus\overline D)}
=O(\tau e^{-\tau\sqrt{\rho/\mu}\,\text{dist}\,(D,B)}+e^{-\tau T}).
$$
This together with the first equation on (2.2) gives
$$\displaystyle
\int_{\Omega\setminus\overline D}\mbox{\boldmath $F$}\cdot\mbox{\boldmath $R$}\,dx
=O(\tau^2 e^{-\tau\sqrt{\rho/\mu}\,\text{dist}\,(D,B)}+\tau e^{-\tau T}).
\tag {3.7}
$$
And also it follows from (2.2), (2.4) and (2.42) with $U=D, \Omega\setminus\overline D$ we obtain
$$
\displaystyle
\int_D\mbox{\boldmath $F$}_0\cdot\mbox{\boldmath $w$}_0dx=O(\tau e^{-\tau\sqrt{\rho/\mu}\,\text{dist}\,(D,B)}+e^{-\tau T})
\tag {3.8}
$$
and
$$
\displaystyle
\int_{\Omega\setminus\overline D}(\mbox{\boldmath $F$}_0-\mbox{\boldmath $F$})
\cdot\mbox{\boldmath $w$}_0\,dx
=O(\tau e^{-\tau\sqrt{\rho/\mu}\,\text{dist}\,(\Omega,B)}+e^{-\tau T}).
\tag {3.9}
$$
Applying (3.7), (3.8) and (3.9) to the right-hand side on (2.10), we obtain
$$
\displaystyle
{\cal R}^1(\tau)
=O(e^{-\tau T}(\tau^2 e^{-\tau\sqrt{\rho/\mu}\,\text{dist}\,(D,B)}+\tau e^{-\tau\sqrt{\rho/\mu}\,\text{dist}\,(\Omega,B)}
+\tau e^{-\tau T})).
\tag {3.10}
$$

Note that ${\cal R}^2(\tau)$ becomes
$$\displaystyle
{\cal R}^2(\tau)=e^{-\tau T}\int_{\Omega\setminus\overline D}h\Sigma\,dx.
\tag {3.11}
$$
From (3.6) we have
$$\displaystyle
e(\tau)
=O(\tau^5 e^{-2\tau\sqrt{\rho/\mu}\,\text{dist}\,(D,B)}+\tau^3 e^{-2\tau T}).
$$
This together with (2.15) yields
$$\displaystyle
\Vert\Sigma\Vert_{L^2(\Omega\setminus\overline D)}
=O(\tau^2e^{-\tau\sqrt{\rho/\mu}\,\text{dist}\,(D,B)}+\tau e^{-\tau T}).
$$
Then, noting $h_0=0$ by (2.4) and $\Xi_0=0$, from the second equation on (2.2) and (3.11), we obtain
$$
\displaystyle
\frac{1}{\theta_0\tau}{\cal R}^2(\tau)
=O(e^{-\tau T}(\tau  e^{-\tau\sqrt{\rho/\mu}\,\text{dist}\,(D,B)}
+e^{-\tau T})).
$$
Now a combination of this and (3.10) gives 
$$\displaystyle
{\cal R}(\tau)
=O(\tau^2e^{-\tau T}e^{-\tau\sqrt{\rho/\mu}\,\text{dist}\,(D,B)}+\tau e^{-\tau T}e^{-\tau\sqrt{\rho/\mu}\,\text{dist}\,(\Omega,B)}
+\tau e^{-2\tau T}).
$$
Finally, applying the trivial estimate $\text{dist}\,(\Omega,B)\le\text{dist}\,(D,B)$ to the first term on 
this right-hand side, we obtain (3.5).

\noindent
$\Box$

\section{Some corollaries}

Some remarks on corollaries of Theorem 1.1 should be mentioned.

We introduce another indicator function defined by the expression
$$\begin{array}{ll}
\displaystyle
I^s(\tau;\mbox{\boldmath $v$}_s,0)
=
\int_{\partial\Omega}
s(\mbox{\boldmath $w$}_{s0},0)\mbox{\boldmath $\nu$}
\cdot
(\mbox{\boldmath $w$}-\mbox{\boldmath $w$}_{s0})\,dS,
&
\tau>0,
\end{array}
$$
where $\mbox{\boldmath $w$}_{s0}\in H^1(\Bbb R^3)^3$ is given by (2.36) and the second equation on (A.3) explicitly.

\proclaim{\noindent Corollary 4.1.}
Theorem 1.1 remains valid if $I^1(\tau;\mbox{\boldmath $v$}_s,0)$ is replaced with 
$I^s(\tau;\mbox{\boldmath $v$}_s,0)$.

\endproclaim

{\it\noindent Proof.}
Write
$$\begin{array}{ll}
\displaystyle
s(\mbox{\boldmath $w$}_{0},0)\mbox{\boldmath $\nu$}
\cdot
(\mbox{\boldmath $w$}-\mbox{\boldmath $w$}_{0})
&
\displaystyle
=s(\mbox{\boldmath $w$}_{s0},0)\mbox{\boldmath $\nu$}
\cdot
(\mbox{\boldmath $w$}-\mbox{\boldmath $w$}_{0})
+
s(\mbox{\boldmath $w$}_{0}-\mbox{\boldmath $w$}_{s0},0)
\mbox{\boldmath $\nu$}
\cdot
(\mbox{\boldmath $w$}-\mbox{\boldmath $w$}_{0})\\
\\
\displaystyle
&
\displaystyle
=
s(\mbox{\boldmath $w$}_{s0},0)\mbox{\boldmath $\nu$}\cdot
(\mbox{\boldmath $w$}-\mbox{\boldmath $w$}_{s0})
+s(\mbox{\boldmath $w$}_{s0},0)\mbox{\boldmath $\nu$}
\cdot(\mbox{\boldmath $w$}_{s0}-\mbox{\boldmath $w$}_{0})
\\
\\
\displaystyle
&
\displaystyle
\,\,\,
+s(\mbox{\boldmath $w$}_0-\mbox{\boldmath $w$}_{s0},0)
\mbox{\boldmath $\nu$}
\cdot
(\mbox{\boldmath $w$}-\mbox{\boldmath $w$}_0)
\\
\\
\displaystyle
&
\displaystyle
=s(\mbox{\boldmath $w$}_{s0},0)\mbox{\boldmath $\nu$}\cdot
(\mbox{\boldmath $w$}-\mbox{\boldmath $w$}_{s0})
\\
\\
\displaystyle
&
\,\,\,
-e^{-\tau T}s(\mbox{\boldmath $w$}_{s0},0)\mbox{\boldmath $\nu$}
\cdot\mbox{\boldmath $\epsilon$}_{s}
+e^{-\tau T}s(\mbox{\boldmath $\epsilon$}_s,0)
\mbox{\boldmath $\nu$}
\cdot
(\mbox{\boldmath $w$}-\mbox{\boldmath $w$}_0),
\end{array}
\tag {4.1}
$$
where $\mbox{\boldmath $\epsilon$}_s=e^{\tau T}(\mbox{\boldmath $w$}_0-\mbox{\boldmath $w$}_{s0})$.
A combination of second estimate on (2.28) and (3.6), we have
$$\displaystyle
\Vert\mbox{\boldmath $w$}-\mbox{\boldmath $w$}_0\Vert_{H^{1/2}(\partial\Omega)}=
O(\tau^2 e^{-\tau\sqrt{\rho/\mu}\,\text{dist}\,(D,B)}
+\tau e^{-\tau T}
).
\tag {4.2}
$$
From (2.46) together with $\Vert\mbox{\boldmath $F$}_0\Vert_{L^2(\Bbb R^3)}=O(\tau)$
and (2.47) we have $\Vert\Delta\mbox{\boldmath $\epsilon$}_s\Vert_{L^2(\Bbb R^3)}=O(\tau)$
and hence 
$$\displaystyle
\Vert\mbox{\boldmath $\epsilon$}_s\Vert_{H^2(\Bbb R^3)}=O(\tau).
\tag {4.3}
$$
This gives
$$
\displaystyle
\left\Vert
s(\mbox{\boldmath $\epsilon$}_s,0)\mbox{\boldmath $\nu$}\right\Vert_{H^{1/2}(\partial\Omega)}=O(\tau).
$$
A combination of this and (4.2) gives
$$\displaystyle
\int_{\partial\Omega}
s(\mbox{\boldmath $\epsilon$}_0,0)\mbox{\boldmath $\nu$}
\cdot(\mbox{\boldmath $w$}-\mbox{\boldmath $w$}_0)\,dS
=O(\tau^3e^{-\tau\sqrt{\rho/\mu}\,\text{dist}\,(D,B)}
+\tau^2e^{-\tau T}
).
$$
And also from (2.30) and (4.3) we have
$$
\displaystyle
\int_{\partial\Omega}s(\mbox{\boldmath $w$}_{s0},0)\mbox{\boldmath $\nu$}
\cdot\mbox{\boldmath $\epsilon$}_s
\,dS
=O(\tau e^{-\tau\sqrt{\rho/\mu}\,\text{dist}\,(\Omega,B)}).
$$
Using $\text{dist}\,(\Omega,B)<\text{dist}\,(D,B)$ and
Integrating both sides on (4.1) and applying these,
we obtain
$$
\displaystyle
I^1(\tau;\mbox{\boldmath $v$}_s,0)
=I^s(\tau;\mbox{\boldmath $v$}_s,0)
+O(\tau^3e^{-\tau(T+\sqrt{\rho/\mu}\,\text{dist}\,(\Omega,B))}+
\tau^2e^{-2\tau T}).
$$
The bound for the last term on this right-hand side is essentially
same as that of ${\cal R}(\tau)$ given in (3.5).
Thus, we can easily see that all the statements of Theorem 1.1 are transplanted.

\noindent
$\Box$

And also one can localize the place where the data are corrected.
Given $M>0$ define the localized indicator function by the formula
$$\begin{array}{ll}
\displaystyle
I^s(\tau;\mbox{\boldmath $v$}_s,0;M)
=
\int_{\partial\Omega(B,M)}s(\mbox{\boldmath $w$}_{s0},0)\mbox{\boldmath $\nu$}
\cdot
(\mbox{\boldmath $w$}-\mbox{\boldmath $w$}_{s0})\,dS,
&
\tau>0,
\end{array}
$$
where
$$\displaystyle
\partial\Omega(B,M)
=\{x\in\partial\Omega\,\vert\, d_B(x)<M\}
$$
and $d_B(x)=\inf_{y\in B}\vert y-x\vert$.

The second corollary of Theorem 1.1 is the following.

\proclaim{\noindent Corollary 4.2.}
Let $M$ satisfy
$$\displaystyle
\text{dist}\,(D,B)<M.
\tag {4.4}
$$
Let $T$ satisfy
$$\displaystyle
T\ge \sqrt{\frac{\rho}{\mu}}\left(2M-\text{dist}\,(\Omega,B)\right).
\tag {4.5}
$$
Then, statement (i) in Theorem 1.1 remains valid if $I^1(\tau;\mbox{\boldmath $v$}_s,0)$ is replaced with
$I^s(\tau;\mbox{\boldmath $v$}_s,0;M)$.

\endproclaim

{\it\noindent Proof.}
From the expression (2.30), we have
$$\displaystyle
\left\Vert
s(\mbox{\boldmath $w$}_{s0},0)\mbox{\boldmath $\nu$}
\right\Vert_{L^2(\partial\Omega\setminus\partial\Omega(M,B))}
=O(\tau e^{-\tau\sqrt{\rho/\mu}\,M}).
$$
It follows from (4.2) and (4.3) that
$$\displaystyle
\Vert\mbox{\boldmath $w$}-\mbox{\boldmath $w$}_{s0}\Vert_{L^2(\partial\Omega)}
=O(\tau^2 e^{-\tau\sqrt{\rho/\mu}\,\text{dist}\,(D,B)}+\tau e^{-\tau T}
).
$$
A combination of these gives
$$\displaystyle
I^s(\tau;\mbox{\boldmath $v$}_0,0)=I^s(\tau;\mbox{\boldmath $v$}_s,0;M)
+O(\tau^3e^{-\tau\sqrt{\rho/\mu}\,\text{dist}\,(D,B)}e^{-\tau\sqrt{\rho/\mu}\, M}
+\tau^2 e^{-\tau T}e^{-\tau\sqrt{\rho/\mu}\, M}
).
$$
Then, we can check the validity of the statement in Corollary 4.2 
by using Corollary 4.1 and the following facts.

$\bullet$  One can write
$$\left\{\begin{array}{l}
\displaystyle
e^{2\tau\sqrt{\rho/\mu}\,\text{dist}\,(D,B)}\tau^3e^{-\tau\sqrt{\rho/\mu}\,\text{dist}\,(D,B)}e^{-\tau\sqrt{\rho/\mu}\,M}
=\tau^3 e^{-\tau\sqrt{\rho/\mu}\,(M-\text{dist}\,(D,B))},\\
\\
\displaystyle
e^{2\tau\sqrt{\rho/\mu}\,\text{dist}\,(D,B)}\tau^2 e^{-\tau T}e^{-\tau\sqrt{\rho/\mu}\,M}
=\tau^2 e^{-\tau\{T+\sqrt{\rho/\mu}\,(M-2\text{dist}\,(D,B))\}}.
\end{array}
\right.
$$

$\bullet$  a combination of (4.4) and  (4.5) implies
$$
\displaystyle
T+\sqrt{\frac{\rho}{\mu}}\,\left(M-2\text{dist}\,(D,B)\right)>0.
$$

\noindent
$\Box$

{\bf\noindent Remark 4.1.}
Corollaries 4.1 and 4.2 remain valid if
the function $\mbox{\boldmath $w$}_{s0}$ in $\mbox{\boldmath $w$}-\mbox{\boldmath $w$}_{s0}$ 
of $I^s(\tau;\mbox{\boldmath $v$}_s,0)$ and
$I^s(\tau;\mbox{\boldmath $v$}_s,0;M)$
is replaced with $\mbox{\boldmath $w$}_0^m(\,\cdot\,,\mbox{\boldmath $v$}_s,0)$ which depends on $T$
and given by the second equation on (1.8) with $(\mbox{\boldmath $v$},\Theta)=(\mbox{\boldmath $v$}_s,0)$.
This is the original style stated in \cite{IE4} for the scalar wave equation.

\section{Further problems}

In Theorem 1.1 we made use of only the first indicator function $I^1(\tau;\mbox{\boldmath $v$}_s,0)$
which employs only a single displacement field observed on the surface of the body.
It is based on the simple fact that the system (1.2) has special solutions $\mbox{\boldmath $v$}=\mbox{\boldmath $v$}_s$
with $\nabla\cdot\mbox{\boldmath $v$}=0$
and $\Theta=0$.  This choice yields that the second indicatior function $I^2(\tau;\mbox{\boldmath $v$}_s,0)$
becomes identially $0$ and thus does not enable us to extract any information about the cavity $D$ from a single temperature field
observed on the surface of the body.  Needless to say there are other choices of special solutions
of the system (1.2).  For example, just solve (1.2) with the initial conditions
$$\displaystyle
\left\{
\begin{array}{ll}
\displaystyle
\mbox{\boldmath $v$}(x,0)=\mbox{\boldmath $0$} & \text{in $\Bbb R^3$,}\\
\\
\displaystyle
\partial_t\mbox{\boldmath $v$}(x,0)=\mbox{\boldmath $0$} & \text{in $\Bbb R^3$,}\\
\\
\displaystyle
\Theta(x,0)=(\eta-\vert x-p\vert)^2\chi_B(x) & \text{in $\Bbb R^3$}
\end{array}
\right.
$$
and take $f$ and $\mbox{\boldmath $G$}$ in (1.1) given by (1.4) and (1.5).
Then it would be interested to clarify the information about the cavity
contained in both the first and second indicator functions $I^1(\tau;\mbox{\boldmath $v$},\Theta)$
and $I^2(\tau;\mbox{\boldmath $v$},\Theta)$.  This needs further analysis 
on the behaviour of the {\it reflected solutions} $\mbox{\boldmath $R$}$ and $\Sigma$
which are solutions of (2.11) and (2.16) and appear in the representation
formula of the indicator functions (2.7) and (2.13).  For the scalar wave equation case the analysis of the effect of the corresponding reflected solution
on the asymptotic behaviour of the indicator function has been done in \cite{ICA, IEO3, IEE, IMP} 
and see also \cite{IMax} for the Maxwell system.  The governing equation of the present problem is a coupled system of the elastic
wave and heat equations and it seems more difficult.  We leave it for future research.
Finally, it should be mentioned that the numerical implementation of the method of this paper belongs to our future plan.

$\quad$

\centerline{{\bf Acknowledgments}}
The author was partially supported by Grant-in-Aid for
Scientific Research (C)(No. 17K05331) of Japan  Society for
the Promotion of Science.  
The author thanks Hiromichi Itou for pointing out the reference \cite{D} 
and useful discussions.
Some part of this work was initiated when the author
stayed in ICUB (The Research Institute of the University of Bucharest), Bucharest, Romania for 1th Nov.-18th Nov. in 2016 with support 
by ICUB the visiting professors fellowship program.

\section{Appendix}

\subsection{Explicit computation of integrals on (2.38) and (2.41)}

The following formulae have been derived in \cite{IE4} as (A.5) and (A.6) therein.
Note that the first formula on (A.1) below is also a result of the mean value theorem for the modified Helmholtz equation \cite{CH}.

\proclaim{\noindent Proposition A.1(\cite{IE4})}
Let $x\in\Bbb R^3\setminus\overline B$.
Then, we have
$$
\left\{
\begin{array}{l}
\displaystyle
I_0(x;\tau)=\frac{4\pi\varphi_0(\tau\eta)}{\tau^3}\frac{e^{-\tau\vert x-p\vert}}{\vert x-p\vert},\\
\\
\displaystyle
I_1(x;\tau)
=\frac{4\pi\varphi_1(\tau\eta)}{\tau^4}
\frac{e^{-\tau\vert x-p\vert}}{\vert x-p\vert},
\end{array}
\right.
\tag {A.1}
$$
where
$$\left\{\begin{array}{l}
\varphi_0(s)=s\cosh s-\sinh s,
\\
\\
\displaystyle
\varphi_1(s)=\left(s^2+2\right)\cosh s
-2s\sinh s
-2.
\end{array}
\right.
\tag {A.2}
$$

\endproclaim

In this appendix, we add the following formulae.

\proclaim{\noindent Proposition A.2.}
We have
$$\left\{
\begin{array}{l}
\displaystyle
\mbox{\boldmath $I$}_0(x;\tau)
=
\pi \frac{e^{-\tau\,\vert x-p\vert}}{\vert x-p\vert^2}
K_{\tau}^0\,(\vert x-p\vert,\tau\eta)
\frac{x-p}{\vert x-p\vert},
\\
\\
\displaystyle
\mbox{\boldmath $I$}_1(x;\tau)
=\frac{4\pi}{\tau^2}
\frac{e^{-\tau\,\vert x-p\vert}}{\vert x-p\vert^2}
K_{\tau}^1(\vert x-p\vert,\tau\eta)\frac{x-p}{\vert x-p\vert},
\\
\\
\displaystyle
I_2(x;\tau)
=\frac{4\pi\varphi_2(\tau\eta)}{\tau^5}\frac{e^{-\tau\vert x-p\vert}}{\vert x-p\vert},
\end{array}
\right.
\tag {A.3}
$$
where
$$\begin{array}{ll}
\displaystyle
K_{\tau}^0(\xi,s)
&
\displaystyle
=\frac{2}{\tau}\left\{
\left(1-\frac{1}{\tau}\right)\xi^2-\frac{2}{\tau^2}\xi-\frac{2}{\tau^3}
\right\}(\cosh s-1)\\
\\
\displaystyle
&
\displaystyle
\,\,\,
+\frac{2}{\tau^3}\left(1-\frac{1}{\tau}\right)
\left\{(s^2+2)\cosh s-2s\sinh s-2\right\}\\
\\
\displaystyle
&
\displaystyle
\,\,\,
+\frac{4}{\tau^3}(s\sinh s-\cosh s),
\end{array}
\tag {A.4}
$$
$$\begin{array}{ll}
\displaystyle
K_{\tau}^1(\xi,s)
&
\displaystyle
=\left\{\left(1-\frac{1}{\tau}\right)\xi^2-\frac{2}{\tau^2}
\left(\xi+\frac{1}{\tau}\right)\right\}(s\cosh s-\sinh s)\\
\\
\displaystyle
&
\displaystyle
\,\,\,
+\frac{1}{\tau^2}\left(1-\frac{1}{\tau}\right)
\left\{s^2(s+6)\cosh s-3(s^2+2)\sinh s\right\}\\
\\
\displaystyle
&
\displaystyle
\,\,\,
+\frac{2}{\tau^2}
\left(\xi+\frac{1}{\tau}\right)
\left\{(s^2+1)\sinh s-2s\cosh s\right\}.
\end{array}
\tag {A.5}
$$
and
$$
\displaystyle
\varphi_2(s)
=s^2(s+6)\cosh s
-3(s^2+2)\sinh s.
\tag {A.6}
$$

\endproclaim

{\it\noindent Proof.}
It suffices to consider the case when $p=0$.
Let $\mbox{\boldmath $A$}_x$ be the orthogonal matrix such that
$(\mbox{\boldmath $A$}_x)^Tx=\vert x\vert\mbox{\boldmath $e$}_3$.
The change of variables $y=r\omega\,(0<r<\eta,\,\omega\in S^2)$ and a rotation give us
$$\begin{array}{ll}
\displaystyle
\mbox{\boldmath $I$}_0(x;\tau) & \displaystyle =\int_0^{\eta}r^2 dr\int_{S^2}\frac{e^{-\tau\vert x-r\omega\vert}}{\vert x-r\omega\vert}
\,\omega d\omega
\\
\\
\displaystyle
& \displaystyle =\int_0^{\eta}r^2 dr
\int_{S^2}
\frac{\displaystyle e^{-\tau\vert \vert x\vert\mbox{\boldmath $e$}_3-r\omega\vert}}
{\displaystyle\vert\vert x\vert\mbox{\boldmath $e$}_3-r\omega\vert}\,\mbox{\boldmath $A$}_x\omega d\omega
\\
\\
\displaystyle
& \displaystyle =
\int_0^{\eta}r^2dr\int_0^{2\pi}d\theta
\int_0^{\pi}\sin\varphi d\varphi
\frac{\displaystyle e^{-\tau\sqrt{\vert x\vert^2-2r\vert x\vert\cos\varphi+r^2}}}
{\displaystyle\sqrt{\vert x\vert^2-2r\vert x\vert\cos\varphi+r^2}}
\mbox{\boldmath $A$}_x
\left(\begin{array}{c}
\sin\varphi\cos\theta\\
\\
\displaystyle
\sin\varphi\sin\theta\\
\\
\displaystyle
\cos\varphi
\end{array}
\right)
\\
\\
\displaystyle
&
\displaystyle
=2\pi
\int_0^{\eta}r^2dr
\int_0^{\pi}\sin\varphi d\varphi
\frac{\displaystyle e^{-\tau\sqrt{\vert x\vert^2-2r\vert x\vert\cos\varphi+r^2}}}
{\displaystyle\sqrt{\vert x\vert^2-2r\vert x\vert\cos\varphi+r^2}}
\cos\varphi\mbox{\boldmath $A$}_x\mbox{\boldmath $e$}_3
\\
\\
\displaystyle
& \displaystyle =2\pi\int_0^{\eta}U(\vert x\vert,r)r^2 dr\frac{x}{\vert x\vert},
\end{array}
\tag {A.7}
$$
where 
$$\begin{array}{lll}
\displaystyle
U(\xi,r)
=\int_0^{\pi}
\frac{\displaystyle e^{-\tau\sqrt{\xi^2-2r\xi\cos\varphi+r^2}}}
{\displaystyle\sqrt{\xi^2-2r\xi\cos\varphi+r^2}}\sin\varphi \cos\varphi\,d\varphi, & \xi>\eta, & 0<r<\eta.
\end{array}
$$
Fix $\xi\in]\eta,\,\infty [$ and $r\in]0,\,\eta[$.
The change of variable
$$\displaystyle
s=\sqrt{\xi^2-2r\xi\cos\varphi+r^2},\,\varphi\in]0,\,\pi[
$$
gives
$$\displaystyle
\sin\varphi\,\cos\varphi d\varphi
=\frac{\xi^2+r^2-s^2}{2r^2\xi^2}s\,ds.
$$
Thus we have
$$\begin{array}{ll}
\displaystyle
U(\xi,r)
& \displaystyle
=\frac{1}{2r^2\xi^2}\int_{\xi-r}^{\xi+r}e^{-\tau s}(\xi^2+r^2-s^2)ds
\\
\\
\displaystyle
& \displaystyle =
\frac{1}{2r^2\xi^2}
\left\{(\xi^2+r^2)(-e^{-\tau(\xi+r)}+e^{-\tau(\xi-r)})
-\int_{\xi-r}^{\xi+r}e^{-\tau s}s^2ds\right\}.
\end{array}
$$
Using the formula
$$\begin{array}{ll}
\displaystyle
\int_{\xi-r}^{\xi+r}e^{-\tau s}s^2ds
&
\displaystyle
=-\frac{1}{\tau}
\left\{(\xi+r)^2+\frac{2}{\tau}(\xi+r)+\frac{2}{\tau^2}\right\}e^{-\tau(\xi+r)}\\
\\
\displaystyle
&
\displaystyle
\,\,\,
+\frac{1}{\tau}
\left\{(\xi-r)^2+\frac{2}{\tau}(\xi-r)+\frac{2}{\tau^2}\right\}e^{-\tau(\xi-r)},
\end{array}
$$
we obtain
$$
\displaystyle
U(\xi,r)
=\frac{e^{-\tau, \xi}}{2\xi^2 r^2}
\left(A(\xi,r)
e^{\tau r}
-
B(\xi,r)
e^{-\tau r}
\right)
$$
and
$$\left\{
\begin{array}{l}
\displaystyle
A(\xi,r)=\left(\xi^2+r^2\right)-\frac{1}{\tau}(\xi-r)^2-\frac{2}{\tau^2}(\xi-r)-\frac{2}{\tau^3},
\\
\\
\displaystyle
B(\xi,r)=\left(\xi^2+r^2\right)-\frac{1}{\tau}(\xi+r)^2-\frac{2}{\tau^2}(\xi+r)-\frac{2}{\tau^3}.
\end{array}
\right.
$$

From this and (A.7) we obtain
$$
\displaystyle
\,\,\,\,\,\,
\mbox{\boldmath $I$}_0(x;\tau)
=
\pi e^{-\tau\xi}
\int_0^{\eta}
\left(
A(\xi,r)e^{\tau r}-
B(\xi,r)
e^{-\tau r}
\right) dr
\vert_{\xi=\vert x\vert}\frac{x}{\vert x\vert^3}.
\tag {A.8}
$$
We have
$$\left\{
\begin{array}{l}
\displaystyle
\int_0^{\eta}e^{\pm\tau r}rdr
=\pm\frac{1}{\tau}\left(\eta\mp\frac{1}{\tau}\right)e^{\pm\tau\eta}+\frac{1}{\tau^2},
\\
\\
\displaystyle
\int_0^{\eta}e^{\pm\tau r}r^2dr
=\pm\frac{1}{\tau}
\left(\eta^2\mp\frac{2\eta}{\tau}+\frac{2}{\tau^2}\right)e^{\pm\tau \eta}
\mp\frac{2}{\tau^3},
\\
\\
\displaystyle
\int_0^{\eta}e^{\pm\tau r}r^3dr
=
\pm\frac{1}{\tau}e^{\pm\tau \eta}\eta^3
-
\frac{3}{\tau^2}
\left(\eta^2\mp\frac{2\eta}{\tau}+\frac{2}{\tau^2}\right)e^{\pm\tau\eta}
+\frac{6}{\tau^4}
\end{array}
\right.
$$
and thus
$$
\left\{
\begin{array}{l}
\displaystyle
\int_0^{\eta}(e^{\tau r}-e^{-\tau r})rdr
=\frac{2}{\tau^2}(\tau\eta\cosh\tau\eta-\sinh\tau\eta),
\\
\\
\displaystyle
\int_0^{\eta}(e^{\tau r}+e^{-\tau r})rdr
=\frac{2}{\tau^2}(\tau\eta\sinh\tau\eta-\cosh\tau\eta),
\\
\\
\displaystyle
\int_0^{\eta}(e^{\tau r}+e^{-\tau r})r^2dr
=
\frac{2}{\tau^3}
\left\{(\eta^2\tau^2+1)\sinh\tau\eta-2\eta\tau\cosh\tau\eta\right\},
\\
\\
\displaystyle
\int_0^{\eta}(e^{\tau r}-e^{-\tau r})r^2dr
=
\frac{2}{\tau^3}
\left\{(\eta^2\tau^2+2)\cosh\tau\eta-2\eta\tau\sinh\tau\eta-2\right\},
\\
\\
\displaystyle
\int_0^{\eta}(e^{\tau r}-e^{-\tau r})r^3dr
=\frac{2}{\tau^4}
\left\{\eta^2\tau^2(\eta\tau+6)\cosh\tau\eta
-3(\eta^2\tau^2+2)\sinh\tau\eta\right\}.
\end{array}
\right.
\tag {A.9}
$$
Now writing
$$\left\{\begin{array}{l}
\displaystyle
A(\xi,r)
=\left(1-\frac{1}{\tau}\right)\xi^2-\frac{2}{\tau^2}\xi-\frac{2}{\tau^3}
+
\left(1-\frac{1}{\tau}\right)r^2
+\frac{2}{\tau}
\left(\xi+\frac{1}{\tau}\right)r,
\\
\\
\displaystyle
B(\xi,r)
=\left(1-\frac{1}{\tau}\right)\xi^2-\frac{2}{\tau^2}\xi-\frac{2}{\tau^3}
+
\left(1-\frac{1}{\tau}\right)r^2
-\frac{2}{\tau}
\left(\xi+\frac{1}{\tau}\right)r,
\end{array}
\right.
\tag {A.10}
$$
from this and (A.9) we obtain
$$\begin{array}{l}
\displaystyle
\,\,\,\,\,\,
\int_0^{\eta}
(A(\xi,r)e^{\tau r}-B(\xi,r)e^{-\tau r}) dr\\
\\
\\
\displaystyle
=
\left\{
\left(1-\frac{1}{\tau}\right)\xi^2-\frac{2}{\tau^2}\xi-\frac{2}{\tau^3}
\right\}
\int_0^{\eta}(e^{\tau r}-e^{-\tau r})\,dr
\\
\\
\displaystyle
\,\,\,
+
\left(1-\frac{1}{\tau}\right)
\int_0^{\eta}(e^{\tau r}-e^{-\tau r})r^2\,dr
\\
\\
\displaystyle
\,\,\,+\frac{2}{\tau}
\int_0^{\eta}(e^{\tau r}+e^{-\tau r})r\,dr\\
\\
\displaystyle
=\frac{2}{\tau}\left\{
\left(1-\frac{1}{\tau}\right)\xi^2-\frac{2}{\tau^2}\xi-\frac{2}{\tau^3}
\right\}(\cosh\tau\eta-1)\\
\\
\displaystyle
+\frac{2}{\tau^3}\left(1-\frac{1}{\tau}\right)
\left\{(\eta^2\tau^2+2)\cosh\tau\eta-2\eta\tau\sinh\tau\eta-2\right\}\\
\\
\displaystyle
+\frac{4}{\tau^3}(\tau\eta\sinh\tau\eta-\cosh\tau\eta)
\end{array}
$$
A combination of this and (A.8) gives the desired formula for $\mbox{\boldmath $I$}_0(x;\tau)$.

Using the same changes of variable as used in the computation of $\mbox{\boldmath $I$}_0(x;\tau)$,
we have
$$
\displaystyle
\,\,\,\,\,\,
\mbox{\boldmath $I$}_1(x;\tau)
=
\pi e^{-\tau\xi}
\int_0^{\eta}
\left(
A(\xi,r)e^{\tau r}-
B(\xi,r)
e^{-\tau r}
\right)r dr
\vert_{\xi=\vert x\vert}\frac{x}{\vert x\vert^3}.
\tag {A.11}
$$
Then, from (A.9) and (A.10) we obtain
$$\begin{array}{l}
\displaystyle
\,\,\,\,\,\,
\int_0^{\eta}
(A(\xi,r)e^{\tau r}-B(\xi,r)e^{-\tau r})r dr\\
\\
\\
\displaystyle
=
\frac{2}{\tau^2}
\left\{
\left(1-\frac{1}{\tau}\right)\xi^2-\frac{2}{\tau^2}\xi-\frac{2}{\tau^3}
\right\}
(\tau\eta\cosh\tau\eta-\sinh\tau\eta)\\
\\
\displaystyle
\,\,\,
+\frac{2}{\tau^4}
\left(1-\frac{1}{\tau}\right)\left\{\eta^2\tau^2(\eta\tau+6)\cosh\tau\eta
-3(\eta^2\tau^2+2)\sinh\tau\eta\right\}\\
\\
\displaystyle
\,\,\,+\frac{4}{\tau^4}
\left(\xi+\frac{1}{\tau}\right)\left\{(\eta^2\tau^2+1)\sinh\tau\eta-2\eta\tau\cosh\tau\eta\right\}.
\end{array}
$$
A combination of this and (A.11) gives the dersired formula for $\mbox{\boldmath $I$}_1(x;\tau)$.

Similarily doing as above, 
one has the expression
$$
\displaystyle
\,\,\,\,\,\,
I_2(x;\tau)
=
\frac{2\pi}{\tau}\frac{e^{-\tau\vert x\vert}}{\vert x\vert}
\int_0^{\eta}
(
e^{\tau r}-
e^{-\tau r})r^3 dr.
$$
Now (A.9) gives the dersired formula for $I_2(x;\tau)$.

\noindent
$\Box$

\subsection{Lower estimates for volume integrals}

Recall that $B$ is an open ball centered at $p$ with radius $\eta$ and satisfies
$\overline B\cap\overline D=\emptyset$: we have $\text{dist}\,(D,B)=d_{\partial D}(p)-\eta$.
In the following lemma it is assumed that $\partial D$ is $C^2$.

\proclaim{\noindent Lemma A.1.}
There exist positive constants $C$, $\tau_0$ and $\kappa$ such that, for all $\tau\ge\tau_0$
$$\displaystyle
\tau^{\kappa}e^{2\tau\sqrt{\rho/\mu}\,\text{dist}\,(D,B)}
\int_D e^{-2\tau\sqrt{\rho/\mu}\,(\vert x-p\vert-\eta)}
\left\vert\frac{x-p}{\vert x-p\vert}\times\mbox{\boldmath $a$}\right\vert^2 dx
\ge C.
\tag {A.12}
$$

\endproclaim

{\it\noindent Proof.}
Choose a point $q\in\partial D$ such that $\vert q-p\vert=d_{\partial D}(p)$.  Since $\partial D$ is $C^2$,
one can find an open ball $B'$ with radius $\delta$ and centered at $q-\delta\mbox{\boldmath $\nu$}_q$ such that
$B'\subset D$ and $\partial B'\cap\partial D=\{q\}$.  Then $\text{dist}\,(B',B)=\text{dist}\,(D,B)$.
Thus, it suffices to prove (A.12) in the case when $D=B'$.

Set $d=d_{\partial D}(p)$ and $-\mbox{\boldmath $\nu$}_q=\mbox{\boldmath $\omega$}_0$.
Let $B''$ be the open ball with radius $d+\delta$ centered at $p$.
First we give a parametrization of the domain $B''\cap B'$.
Given $s\in\,]0,\,\delta[$ we find the set of all unit vectors $\mbox{\boldmath $\omega$}$ satisfying
$p+(d+s)\mbox{\boldmath $\omega$}\in B'$. 
Since the center of $B'$ has the expression $p+(d+\delta)\mbox{\boldmath $\omega$}_{0}$,
this condition is equivalent to the equation
$$\displaystyle
\vert (d+s)\mbox{\boldmath $\omega$}-(d+\delta)\mbox{\boldmath $\omega$}_0\vert<\delta.
$$
This is equivalent to the condition
$$\displaystyle
\mbox{\boldmath $\omega$}\cdot\mbox{\boldmath $\omega$}_0>\frac{(d+\delta)^2+(d+s)^2-\delta^2}{2(d+s)(d+\delta)}.
$$
Thus, we have
$$\displaystyle
B''\cap B'
=\cup_{0<s<\delta}
\left\{p+(d+s)\omega\,\vert\,\omega\in S(s)\right\},
$$
where
$$\displaystyle
S(s)=\left\{\mbox{\boldmath $\omega$}\in S^2\,\vert\,\mbox{\boldmath $\omega$}\cdot\mbox{\boldmath $\omega$}_0
>\frac{(d+\delta)^2+(d+s)^2-\delta^2}{2(d+s)(d+\delta)}\right\}.
$$
Choose two linearly independent vectors $\mbox{\boldmath $b$}$ and $\mbox{\boldmath $c$}$ in such a way that
$\mbox{\boldmath $b$}\cdot\mbox{\boldmath $c$}=0$, $\mbox{\boldmath $b$}\times\mbox{\boldmath $c$}=
\mbox{\boldmath $\omega$}_0$.
We denote by $\theta(s)\in\,]0, \pi/2[$ the unique solution of
$$\displaystyle
\cos\theta
=\frac{(d+\delta)^2+(d+s)^2-\delta^2}{2(d+s)(d+\delta)}.
\tag {A.13}
$$
Given $s\in\,]0,\,\delta[$, $r\in\,]0,\,(d+s)\sin\theta(s)[$ and $\gamma\in\,[0,\,2\pi[$ set
$$\displaystyle
\mbox{\boldmath $\Upsilon$}(s,r,\gamma)
=p+(d+s)\cos\theta(s)\,\mbox{\boldmath $\omega$}_0
+r(\cos\gamma\,\mbox{\boldmath $b$}+\sin\gamma\,\mbox{\boldmath $c$})
+h\,\mbox{\boldmath $\omega$}_0,
$$
where $h$ is an unknown parameter to be determined by the equation
$$\displaystyle
\vert\mbox{\boldmath $\Upsilon$}(s,r,\gamma)-p\vert=d+s.
$$
By solving this, we obtain
$$\displaystyle
h=-(d+s)\cos\theta(s)+\sqrt{(d+s)^2-r^2}.
$$
Thus, we have
$$\displaystyle
\mbox{\boldmath $\Upsilon$}(s,r,\gamma)
=p+\sqrt{(d+s)^2-r^2}\,\mbox{\boldmath $\omega$}_0
+r(\cos\gamma\,\mbox{\boldmath $b$}+\sin\gamma\,\mbox{\boldmath $c$}).
$$
Note that $\vert\mbox{\boldmath $\Upsilon$}(s,r,\gamma)-p\vert=d+s$ and
the unit vector $(\mbox{\boldmath $\Upsilon$}(s,r,\gamma)-p)/(d+s)$
belongs to $S(s)$ for each fixed $s$.  It is easy to check also that the map
$$\displaystyle
G\ni (s,r,\gamma)\longmapsto \mbox{\boldmath $\Upsilon$}(s,r,\gamma)\in B''\cap B'
$$
is bijective, where
$$\displaystyle
G=\left\{(s,r,\gamma)\,\vert\,(s,r)\in G', \gamma\in\,[0,\,2\pi[\,\right\}
$$
and
$$\displaystyle
G'=\left\{(s,r)\,\vert 0<s<\delta,\, 0<r<(d+s)\sin\theta(s)\right\}.
$$
A simple computation gives
$$\displaystyle
\text{det}\,\mbox{\boldmath $\Upsilon$}'(s,r,\gamma)
=\frac{r(d+s)}
{\sqrt{(d+s)^2-r^2}}.
$$
Here we have
$$\begin{array}{ll}
\displaystyle
(\mbox{\boldmath $\Upsilon$}(s,r,\gamma)-p)\times\mbox{\boldmath $a$}
&
\displaystyle
=\sqrt{(d+s)^2-r^2}\mbox{\boldmath $\omega$}_0\times
\mbox{\boldmath $a$}
+r(\cos\gamma\,\mbox{\boldmath $b$}\times\mbox{\boldmath $a$}
+\sin\gamma\,\mbox{\boldmath $c$}\times\mbox{\boldmath $a$})\\
\\
\displaystyle
&
\displaystyle
\equiv \mbox{\boldmath $A$}(s,r,\gamma)\mbox{\boldmath $a$}.
\end{array}
\tag {A.14}
$$
Using the change of variables formula and $B''\cap B'\subset B'$, we obtain
$$\begin{array}{l}
\displaystyle
\,\,\,\,\,\,
e^{2\tau\sqrt{\rho/\mu}\,\,(d-\eta)}\int_{B'} e^{-2\tau\sqrt{\rho/\mu}\,(\vert x-p\vert-\eta)}
\left\vert\frac{x-p}{\vert x-p\vert}\times\mbox{\boldmath $a$}\right\vert^2 dx
\\
\\
\displaystyle
\ge
\int_0^{\delta}ds\int_0^{(d+s)\sin\theta(s)}dr\int_0^{2\pi}d\gamma
\frac{e^{-2\tau\sqrt{\rho/\mu}\,s}}{d+s}
\left\vert\mbox{\boldmath $A$}(s,r,\gamma)\mbox{\boldmath $a$}\right\vert^2
\frac{r}
{\sqrt{(d+s)^2-r^2}}.
\end{array}
\tag {A.15}
$$
It is easy to see that from (A.14) one gets
$$\begin{array}{l}
\displaystyle
\,\,\,\,\,\,
\int_0^{2\pi}\left\vert\mbox{\boldmath $A$}(s,r,\gamma)\mbox{\boldmath $a$}\right\vert^2 d\gamma
\\
\\
\displaystyle
=2\pi\left\{(d+s)^2-r^2\right\}
\vert\mbox{\boldmath $\omega$}_0\times
\mbox{\boldmath $a$}\vert^2
+\pi r^2
(\vert\mbox{\boldmath $b$}\times\mbox{\boldmath $a$}\vert^2
+\vert\mbox{\boldmath $c$}\times\mbox{\boldmath $a$}\vert^2).
\end{array}
$$
And also
$$
\left\{
\begin{array}{l}
\displaystyle
\int_0^{(d+s)\sin\theta(s)}r\sqrt{(d+s)^2-r^2}dr
=\frac{(d+s)^3}{3}(1-\cos^3\theta(s)),\\
\\
\displaystyle
\int_0^{(d+s)\sin\theta(s)}
\frac{r^3dr}
{\sqrt{(d+s)^2-r^2}}
=(d+s)^3
\left\{(1-\cos\theta(s))
+\frac{1-\cos^3\theta(s)}{3}
\right\}.
\end{array}
\right.
$$
Thus
$$\begin{array}{l}
\displaystyle
\,\,\,\,\,\,
\int_0^{(d+s)\sin\theta(s)}\frac{r dr}{\sqrt{(d+s)^2-r^2}}
\int_0^{2\pi}\left\vert\mbox{\boldmath $A$}(s,r,\gamma)\mbox{\boldmath $a$}\right\vert^2 d\gamma\\
\\
\displaystyle
=\frac{2\pi}{3}(d+s)^3(1-\cos^3\theta(s))\vert\mbox{\boldmath $\omega$}_0\times
\mbox{\boldmath $a$}\vert^2\\
\\
\displaystyle
\,\,\,
+\pi 
(d+s)^3\left\{(1-\cos\theta(s))
-\frac{1-\cos^3\theta(s)}{3}
\right\}
(\vert\mbox{\boldmath $b$}\times\mbox{\boldmath $a$}\vert^2
+\vert\mbox{\boldmath $c$}\times\mbox{\boldmath $a$}\vert^2).
\end{array}
$$
From (A.13) we have
$$\begin{array}{ll}
\displaystyle
(d+s)(1-\cos\theta(s))
&
\displaystyle
=\frac{s(2\delta-s)}{2(d+\delta)}
\\
\\
\displaystyle
&
\displaystyle
\ge\frac{\delta s}{2(d+\delta)}.
\end{array}
$$
This gives
$$
\begin{array}{ll}
\displaystyle
(d+s)^3(1-\cos^3\theta(s))
&
\displaystyle
\ge
\frac{\delta (d+s)^2 s}{2(d+\delta)}
(1+\cos\theta(s)+\cos^2\theta(s))
\\
\\
\displaystyle
\,\,\,
&
\displaystyle
\ge\frac{\delta d(d+s)s}{2(d+\delta)}.
\end{array}
$$
And also
$$\begin{array}{ll}
\displaystyle
(d+s)^3\left\{(1-\cos\theta(s))
-\frac{1-\cos^3\theta(s)}{3}
\right\}
&
\displaystyle
\ge
\frac{\delta (d+s)^2s}{2(d+\delta)}
\left(1-\frac{1+\cos\theta(s)+\cos^2\theta(s)}{3}\right)
\\
\\
\displaystyle
&
\displaystyle
=\frac{\delta (d+s)^2s}{6(d+\delta)}
(1-\cos\theta(s))(2+\cos\theta(s))
\\
\\
\displaystyle
&
\displaystyle
\ge
\frac{\delta^2(d+s)s^2}{12(d+\delta)^2}(2+\cos\theta(s))
\\
\\
\displaystyle
&
\displaystyle
\ge
\frac{\delta^2(d+s)s^2}{6(d+\delta)^2}.
\end{array}
$$
Thus we have
$$\begin{array}{l}
\displaystyle
\,\,\,\,\,\,
\int_0^sds\frac{e^{-2\tau\sqrt{\rho/\mu}\,s}}{d+s}
\int_0^{(d+s)\sin\theta(s)}\frac{r dr}{\sqrt{(d+s)^2-r^2}}
\int_0^{2\pi}\left\vert\mbox{\boldmath $A$}(s,r,\gamma)\mbox{\boldmath $a$}\right\vert^2 d\gamma\\
\\
\displaystyle
\ge
\frac{\pi}{3}\frac{\delta d}{d+\delta}
\int_0^{\delta}se^{-2\tau\sqrt{\rho/\mu}\,s}\,ds
\,\vert\mbox{\boldmath $\omega$}_0\times
\mbox{\boldmath $a$}\vert^2\\
\\
\displaystyle
\,\,\,
+\frac{\pi}{6}
\left(\frac{\delta}{d+\delta}\right)^2
\int_0^{\delta}s^2 e^{-2\tau\sqrt{\rho/\mu}\,s}\,ds
\,
(\vert\mbox{\boldmath $b$}\times\mbox{\boldmath $a$}\vert^2
+\vert\mbox{\boldmath $c$}\times\mbox{\boldmath $a$}\vert^2).
\end{array}
\tag {A.16}
$$
Here we have, for $j=1,2$
$$
\displaystyle
\int_0^{\delta}s^je^{-2\tau\sqrt{\rho/\mu}\,s}\,ds
=\frac{1}{(2\tau\sqrt{\rho/\mu})^{j+1}}+O(\tau^{-1}e^{-2\tau\sqrt{\rho/\mu}\,\delta}).
$$
Therefore, from these and (A.15) and (A.16) one can conclude that:
if $\mbox{\boldmath $\omega$}_0\times\mbox{\boldmath $a$}\not=\mbox{\boldmath $0$}$, then, one can choose $\kappa=2$ in (A.12);
if $\mbox{\boldmath $\omega$}_0\times\mbox{\boldmath $a$}=\mbox{\boldmath $0$}$, then $\mbox{\boldmath $a$}=\pm\mbox{\boldmath $\omega$}_0$
and thus $\vert\mbox{\boldmath $b$}\times\mbox{\boldmath $a$}\vert^2
+\vert\mbox{\boldmath $c$}\times\mbox{\boldmath $a$}\vert^2>0$, and one can choose $\kappa=3$ in (A.12).

\noindent
$\Box$

\proclaim{\noindent Lemma A.2.}
There exist positive constants $C$ and $\tau_0$ such that, for all $\tau\ge\tau_0$
$$\displaystyle
\tau^2e^{2\tau\,\text{dist}\,(D,B)}
\int_De^{-2\tau\,\,(\vert x-p\vert-\eta)}\,dx
\ge C.
$$

\endproclaim

This is proved as follows.
From the proof of Lemma A.1
we have
$$\begin{array}{l}
\displaystyle
\,\,\,\,\,\,
e^{2\tau\,(d-\eta)}\int_{D}e^{-2\tau\,(\vert x-p\vert-\eta)}\,dx\\
\\
\displaystyle
\ge
e^{2\tau\,(d-\eta)}\int_{B'}e^{-2\tau\,(\vert x-p\vert-\eta)}\,dx\\
\\
\displaystyle
\ge
2\pi\int_0^{\delta}\,ds\int_0^{(d+s)\sin\theta(s)}\,dr
e^{-2\tau s}\frac{r(d+s)}{\sqrt{(d+s)^2-r^2}}\\
\\
\displaystyle
=2\pi\int_0^{\delta}(d+s)e^{-2\tau s}\,ds
\int_0^{(d+s)\sin\theta(s)}
\frac{r}{\sqrt{(d+s)^2-r^2}}\,dr\\
\\
\displaystyle
=\pi\int_0^{\delta}(d+s)^2(1-\cos\theta(s))e^{-2\tau s}\,ds\\
\\
\displaystyle
\ge
\pi\int_0^{\delta}(d+s)\cdot\frac{\delta s}{2(d+\delta)}\cdot e^{-2\tau s}\,ds
\\
\\
\displaystyle
\ge
\frac{\pi\delta d}{2(d+\delta)}
\int_0^{\delta}se^{-2\tau s}\,ds.
\end{array}
$$
Since
$$\displaystyle
\int_0^{\delta}se^{-2\tau s}\,ds
=\frac{1}{4\tau^2}\int_0^{2\tau\delta}\xi e^{-\xi}\,d\xi,
$$
we obtain the desired conclusion.

\vskip1cm
\noindent
e-mail address

Masaru Ikehata: ikehata@hiroshima-u.ac.jp

\end{document}